\pdfoutput=1
\documentclass[12pt]{article}
\usepackage{epsfig,latexsym,amsmath,amssymb,bm,cite}
\usepackage{graphicx}
\usepackage{graphicx,picture}
\usepackage{pdfpages}

\usepackage{color}
\usepackage{empheq}
\usepackage{url}
\newcommand*\widefbox[1]{\fbox{\hspace{2em}#1\hspace{2em}}}

\usepackage{hyperref}
\usepackage[bottom]{footmisc}

\newcommand{\be}{\begin{equation}}
\newcommand{\ee}{\end{equation}}
\newcommand{\bea}{\begin{eqnarray}}
\newcommand{\eea}{\end{eqnarray}}
\newcommand{\bem}{\begin{multline}}
\newcommand{\eem}{\end{multline}}
\newcommand{\beg}{\begin{gather}}
\newcommand{\eeg}{\end{gather}}

\newcommand{\ben}{\begin{eqnarray*}}
\newcommand{\een}{\end{eqnarray*}}

\begin{document}

\title{{\bf  A story of balls, randomness and PDEs \\[1.5cm] }}
\author{
\vspace{-1.2in}
{\bf Anastasios Taliotis\thanks{e-mail: taliotis.a@gmail.com}
\thanks{Address: JPMorgan Chase, 25 Bank Street, Canary Wharf, London E14 5JP, UK.}
\thanks{Disclaimer: \hspace{-0.04cm}The current \hspace{-0.04cm}paper \hspace{-0.04cm}consists of a personal \hspace{-0.04cm}work curried out by A. \hspace{-0.04cm}Taliotis, \hspace{-0.04cm}and it expresses \hspace{-0.02cm}his \hspace{-0.02cm}own personal views on the particular topic. In particular, this work is not, in any way, endorsed by or related with JPMorgan Chase (the firm) or its employees or stakeholders or any interests or entities represented by the firm.}
\thanks{Keywords: Probability, recursions, generating functional, PDEs, Laurent expansion, residues, Python code.}}
}

\providecommand{\keywords}[1]{\textbf{\textit{Keywords:}} #1}

\vspace{-3in}



\maketitle

\thispagestyle{empty}

\begin{abstract}\hspace{-0.6cm}
Several differential equations usually appearing in mathematical physics are solved through a power series expansion using the Frobenius method, which reduces in solving difference equations (recursions). 
In this paper, a probability problem is presented whose solution follows a completely reversed but systematic approach. Hence, this work is about illustrating how complex probability problems involving difference equations could be tackled with the more powerful techniques of a better studied and well understood field, that of differential equations. The problem is defined as follows: Inside a box containing $r$ red and $w$ white balls \hspace{-0.02cm}random\hspace{-0.02cm} removals \hspace{-0.04cm}occur. \hspace{-0.04cm}The balls are removed \hspace{-0.04cm}successively \hspace{-0.04cm}according to

\hspace{-0.55cm}the three following rules. Rule I: If a white ball is chosen it is immediately discarded from the box. If a red ball is chosen, it is placed back into the box and a new ball is randomly chosen. The second ball is then removed irrespective of the color. Rule II: Once one ball is removed, the game continues from Rule I. Rule III: The game ends once all the red balls are removed. The first question posed is the determination of the probability that $k$ white balls remain where $k=0,\, 1,\, 2,\, ...,\, w$. Ending the game once all the white balls are removed, a second question is the determination of the probability that k red balls remain where $k=0,\, 1,\, 2,\, ...,\, r$. While inductive solutions are possible, the\hspace{-0.04cm} current approach demonstrates a different
%
%
and algorithmic route. In particular, the law of total probability yields a recursive equation that  
%
is transformed into a linear PDE in two dimensions with inhomogeneous source terms, and suitable boundary conditions that depend on $k$. The PDE solutions, which are found analytically, provide up to a known rescaling factor, the (two) generating functionals of the required probabilities as a function of $r$, $w$ and $k$. Using then the derived functionals, the required probability formulas for any $r$, $w$ and $k$ are finally obtained in a closed form; the probability distributions turn out to be linear combinations of hypergeometric functions of type $_3F_2$.
 Reproducing existing results of the academic literature, which are special but less involved cases of the current completely generic solution, this method is quite generic and adaptable to a large class of problems.

 \end{abstract}
 
 
\thispagestyle{empty}


\newpage

\tableofcontents

\setcounter{page}{1}

\section{Introduction}
During interviews for quantitative researchers in investment banks, hedge funds, asset management firms, fin-tech companies etc, the candidates are often requested to solve various probability problems. The current problem is motivated by such interviews and is a more complicated version of these kind of questions.

\vspace{0.2cm}
\hspace{-0.6cm}{\bf Problem 1}

\vspace{0.2cm}
\hspace{-0.7cm}
{\it Inside a box containing $r$ red and $w$ white balls random removals occur. The balls are removed one by one according to the three following rules. Rule I: If a white ball is chosen it is immediately discarded from the box. If a red ball is chosen, it is placed back into the box and a new ball is randomly chosen. The second ball is then removed irrespective of the color. Rule II: Once a ball is removed, the game continues from Rule I. Rule III: The game ends once all the red balls are removed. The question posed is the determination of the probability that $k$ white balls remain where $k=0,\, 1,\, 2,\, ...,\, w$. Changing rule Rule III into Rule IV: The game ends once all the white balls are removed. The second question posed is the determination of the probability that $k$ red balls remain where $k=0,\, 1,\, 2,\, ...,\, r$.}
\vspace{0.2cm}

This paper in not only about providing the solution to Problem 1, as, according to appendix \ref{Gminus}, a combinatorial-inductive solution\footnote{We are particularly grateful to D. Christofides for providing us with the solution through such an approach.} is (also) possible. This paper is rather about illustrating how difficult problems involving difference equations usually appearing in probability questions, could be tackled systematically with the more powerful techniques of a better studied and well understood field, that of differential equations. As it is known, differential equations appearing in mathematical physics, such as Bessel functions, Hypergeometric functions, Hermite polynomials, to name a few, are solved perturbatively using the Frobenius method. During the course of the solution, recursive relations that connect higher order coefficients with lower order ones appear; that is is difference equations. 

In this problem, we work in the reversed order following an algorithmic recipe explained below and summarized in five steps in section \ref{diss}. That is starting from a probability problem, which is reduced to a linear difference equation in two variables, we deduce an inhomogeneous linear partial differential equation (PDE) with suitable boundary conditions that we solve analytically. The solution, up to a known constant, provides the probability generating functional of the problem at hand. 
Expanding then the functional into a Laurent series, we explicitly obtain the probability solutions in closed form, which, as we will compute through our step by step systematic approach, are given by equations (\ref{finallypk}) and (\ref{finallypk4}) for the Rule III and the Rule IV respectively, and are reproduced below. In particular, the probability $p^{(k)}(r,w)$, where $k$ is the number of the white or red remaining balls, and $r$ and $w$ are the red and the white initial balls inside the box respectively, is given by

\begin{equation*}
\boxed{
\begin{aligned}
p_{III}^{(k)}(r,w) &=k!\frac{ r!(r+w+1)}{(r+k+1)!}\,\, \frac{rw!(r+w-k-1)!}{(r+w)!(w-k)!} \mbox{, Rule III,  }\\ 
p_{IV}^{(k)}(r,w) &=(2k+1)\frac{ r!(r+w+k)!}{(r+k+1)!(w+r)!}\,\, \frac{wr!(r+w-k-1)!}{(r+w)!(r-k)!}\mbox{, Rule IV. }
\end{aligned}
}
\end{equation*}

Satisfying the normalization conditions

\begin{equation*}
\boxed{
\begin{aligned}
\sum_{k=0}^wp_{III}^{(k)}(r,w) =1,\,\,\, \sum_{k=0}^rp_{IV}^{(k)}(r,w) =1,
\end{aligned}
}
\end{equation*}
the two probability formulas consist of the main result of this paper.

The proposed approach combines several beautiful branches of mathematics: probability theory, differential equations, special functions and peripherally, through analytic continuation, residues and Laurent expansions of functions, complex analysis. Most of the steps between equations are either done explicitly or explained in detail. Therefore, being self contained, this paper can be served as an introductory set of lecture notes for advanced undergraduate and graduate students with mathematical background that are interested in applications of the aforementioned fields of mathematics.

For a quick but yet a detailed exposition of the reader to the main ideas of this work, it suffices to restrict to sections \ref{2.1}, \ref{pdeRes} and \ref{special}.
 
\vspace{0.1in}
This work is organized as follows.
\vspace{0.1in}

In section \ref{toy}, using $z$-transformations in two dimensions, we solve an easier problem, Problem 2, that involves algebraic but not differential equations. This serves as an introduction to Problem 1 and to the subsequent sections, and clarifies the formalism and the notation.

\vspace{0.1in}

In section \ref{pde} we return to the original problem where using the law of total probability and $z$-transformations, we eventually obtain the required PDE (of Problem 1). The boundary conditions are discussed in the subsequent sections.

\vspace{0.1in}

Section \ref{special} deals with Problem 3, which is a special case of Problem 1, and whose result is known in the community \cite{MuL}. Hence, we cross-check our method against known cases. 
\vspace{0.1in}

Returning to Problem 1, Rule III, section \ref{genPDE} discusses the boundary conditions and derives the generic solution of the PDE for any remaining white balls $k=0,\,1,\,2\,...,\,w$.

\vspace{0.1in}
Using the solutions of the PDE from section \ref{genPDE}, the first part of section \ref{result} derives the explicit probability formulas for the cases $k=0,\,1,\,2$ (Rule III). Given these solutions and by observing the pattern (there exists an even more systematic approach), the remaining section derives the generic probabilities $p_{III}^{(k)}(r,w) $ as a function of the initial balls $r$ and $w$ and of the remaining white balls $k$, which consists the one of the two main results of this paper. The probabilities $p_{III}^{(k)}(r,w) $ are subsequently being investigated.
\vspace{0.1in}

Armed with the experience from earlier sections and in particular from section \ref{result}, section \ref{result4} outlines and eventually provides the second main result of this paper, namely the formula of $p_{IV}^{(k)}(r,w) $ for  Problem 1, Rule IV. This section is much shorter than sections \ref{genPDE} and \ref{result}.

\vspace{0.1in}

The last section, section \ref{diss}, concludes by discussing the advantages of our approach.

\vspace{0.1in}

Appendices \ref{A}, \ref{B}, \ref{C} and \ref{D} are reserved for longer calculations and for proving several useful identities. Appendix \ref{Z} discusses alternative boundary conditions for Problem 1, Rule III. 
Appendix \ref{Gminus} provides the solution to Problem 1, Rule III through a combinatorial-inducive approach. On the other hand, appendix \ref{E} documents a different and hard question posed in \cite{OaPe}, which serves as another reference for these kind of problems. Lastly, appendix \ref{F} provides a python code that simulates Problem 1 and provides the results for a few cases contrasting them against the two analytically obtained formulas. 

\vspace{0.1in}
\underline{Notation}:
\vspace{0.1in}

We use lower case latin letters for the remaining balls inside the box. In particular, we use the letter $r$ for the red balls and the letter $w$ for the white balls that the box contains at the beginning of the games. The letter $k$ is used for the remaining white or red balls (for Rule III and for Rule IV respectively) after the game ends. We use other lower case latin letters except from $k$, $r$ and $w$ in order to sum over other indices that may appear along the way. Finally, during the $z-$transformations, the letters $z$ and $u$ are associated with the red and the white balls respectively.

We define the Kronecker delta by
\begin{align}
\delta_j^i = \mbox{$1$ if $i=j$ $(i,\,j\,=\,\mathbb{Z})$ and $0$ otherwise}.
\end{align} 

The following definitions apply for Problem 1, Rule III (and for Problem 2).

\begin{subequations} \label{pa2e}
\begin{align}
 \label{pa}
p_{III}^{(k)}(r,w) \equiv & \mbox{ The probability to remain with $k$ white balls} \\ \nonumber
&\mbox{ starting with r red and w white balls} \\ \nonumber
\equiv &\lim_{r'\to\ r} \left(\lim_{w'\to\ w}\left(\lim_{k'\to\ k} \left(p_{III}^{(k')}(r',w')  \right)\right)\right), \\
\label{pb}
p_{III}^{(k)}(0,w) \equiv &\,\delta_w^k,\\ 
\label{pc}
p_{III}^{(0)}(0,0) \equiv &1\,,\\
\label{pbb}
p_{III}^{(k)}(r\geq0,0) \equiv & \, \delta_0^k,\\
\label{pe}
p_{III}^{(k)}(r,w<k) \equiv & \, 0. 
\end{align}
\end{subequations}
Equation (\ref{pc}) is a special case of (\ref{pbb}) but we make it explicit. What it should be emphasized is that equations (\ref{pb})-(\ref{pe}) should be consistent with equation (\ref{pa}); this can be checked explicitly once the formula for $p_{III}^{(k)}(r,w)$ is obtained. In particular, as we will see in section \ref{genPDE}, either equation (\ref{pb}) or equation (\ref{pbb}), under particular modifications that we will derive in what follows, should be used as the boundary conditions for the continuum version of the problem.

The motivation behind the specification of the ordering of the limits in equation (\ref{pa}) is because we would like to write compact expressions for the final probability formulas by giving a meaning to subtle expressions such as $w/(w+r)$ or $r/(w+r)$  when $w=r=0$. Another example would be $(r+w-3)!/(r-1)!$ for $r=0$ and $w=2$ and so forth. The particular ordering of the limits implies that for the base case $r=0$ and $w=0$ (the game starts without any red balls, i.e. the game ends before it begins), the probability  $p_{III}^{(k)}(r=0,w=0)$ is, by our limit convention, equal to one when $k=0$ and zero otherwise. To rephrase, the ordering of the limits essentially defines equations (\ref{pb}) and (\ref{pbb}) and hence, according to the discussion of the previous paragraph, such an ordering defines the boundary conditions of the problem, which are necessary for the uniqueness of the solution. We will see explicitly the necessity of the orderings of the limits and the applications of (\ref{pb})-(\ref{pe}) in what follows.


For Problem 1, Rule IV, the roles of $r$ and $w$ are interchanged. In particular, we have

\begin{subequations} \label{pa2e4}
\begin{align}
 \label{pa4}
p_{IV}^{(k)}(r,w) \equiv & \mbox{ The probability to remain with $k$ red balls} \\ \nonumber
&\mbox{ starting with r red and w white balls} \\ \nonumber
\equiv &\lim_{w'\to\ w} \left(\lim_{r'\to\ r}\left(\lim_{k'\to\ k} \left(p_{IV}^{(k')}(r',w')  \right)\right)\right), \\
\label{pb4}
p_{IV}^{(k)}(r,0) \equiv &\,\delta_r^k,\\ 
\label{pc4}
p_{IV}^{(0)}(0,0) \equiv &1\,,\\
\label{pbb4}
p_{IV}^{(k)}(0,w\geq0) \equiv & \, \delta_0^k,\\
\label{pe4}
p_{IV}^{(k)}(r<k,w) \equiv & \, 0. 
\end{align}
\end{subequations}

Analogous comments apply here as the comments below equation (\ref{pa2e}).

\section{A toy model}\label{toy}

The following problem is given.

\vspace{0.2cm}
\hspace{-0.6cm}
{\bf Problem 2}

\vspace{0.2cm}
\hspace{-0.7cm}
{\it Inside a box containing $r$ red and $w$ white balls random removals occur. The balls are removed randomly one by one. The game ends once all the red balls are removed. The question posed is the determination of the probability that $k$ white balls remain where $k=0,\, 1,\, 2,\, ..., w$.}
\vspace{0.2cm}

\subsection{Solution through $z$-transformations}\label{2.1}

\vspace{0.2cm}

Using the law of total probability by conditioning on the first removal, the following recursive relation is obtained 
\begin{align}\label{rec1}
p^{(k)}(r,w) = \frac{r}{r+w}p^{(k)}(r-1,w) + \frac{w}{r+w}p^{(k)}(r,w-1).
\end{align}

We now make the educated ansatz
\begin{align}\label{ans1}
p^{(k)}(r,w) =  \frac{r!w!}{(r+w)!}f^{(k)}(r,w)
\end{align}
obtaining the much simpler equation for $f^{(k)}$
\begin{align}\label{rec2}
f^{(k)}(r,w) =  f^{(k)}(r-1,w) + f^{(k)}(r,w-1).
\end{align}
This ansatz is motivated because we want to factor out $\frac{r!w!}{(r+w)!}$, which provides the probability to remain with $k=w$ white balls; that is if we keep removing only red balls (equivalently only white balls).

In the next step, we multiply both sides of (\ref{rec2}) by $z^{-r}$ and $u^{-w}$ and sum over $r$ and $w$ from one to infinity. Shifting some indices around and adding and subtracting suitable terms, we eventually obtain
\begin{align}\label{sum1}
&\sum_{r=0}^{\infty} \sum_{w=0} ^{\infty }f^{(k)}_{r,w} z^{-r} u^{-w} -  \sum_{w=0} ^{\infty }f^{(k)}_{0,w} u^{-w} -  \sum_{r=0} ^{\infty }f^{(k)}_{r,0} z^{-r} + f_{0,0} =\\ \nonumber
& \frac{1}{z} \left(  \sum_{r=0}^{\infty} \sum_{w=0} ^{\infty }f^{(k)}_{r,w} z^{-r} u^{-w} - \sum_{r=0} ^{\infty }f^{(k)}_{r,0} z^{-r}  \right) +
\frac{1}{u} \left(  \sum_{r=0}^{\infty} \sum_{w=0} ^{\infty }f^{(k)}_{r,w} z^{-r} u^{-w} - \sum_{w=0} ^{\infty }f^{(k)}_{0,w} u^{-w}  \right)
\end{align}
where we simplified the notation on $f^{(k)}$ setting $f^{(k)}(r,w)=f^{(k)}_{r,w}$.

Defining the generating functional $Y^{(k)}$ by
\begin{align}\label{Yfdef}
Y^{(k)}(z,u) = \sum_{r=0}^{\infty} \sum_{w=0} ^{\infty }f^{(k)}_{r,w} z^{-r} u^{-w} 
\end{align}
equation (\ref{sum1}) yields
\begin{align}\label{Y1}
Y^{(k)}(z,u) = \frac{1}{zu-z-u} \left( u(z-1)  \sum_{r=0} ^{\infty }f^{(k)}_{r,0} z^{-r}  + z(u-1)  \sum_{w=0} ^{\infty }f^{(k)}_{0,w} u^{-w} - z uf^{(k)}_{0,0} \right).
\end{align}

\vspace{0.2cm}
\underline{ Case $k>0$}
\vspace{0.2cm}

In this case, most of the terms in the right hand side of (\ref{Y1}) are zero. In particular, $f^{(k>0)}_{r,0}=0 \, \forall \,r \geq 0$ because the probability to remain with a positive number of white balls $k$ starting without any white balls at all must be zero. Another way to see it is through equation (\ref{pbb}). Also in the view of (\ref{pb}) all $f^{(k)}_{0,w}$'s are zero except from the term $f^{(k)}_{0,k}$, which is equal to one. Hence,
\begin{align}\label{Y1k}
Y^{(k>0)}(z,u) = \frac{z(u-1)}{zu-z-u} \times u^{-k} \mbox{ (because $f^{(k>0)}_{r,0}=0\, \,\forall \,r $ and $f^{(k)}_{0,w}=\delta^k_w$).}
\end{align}
The $f^{(k)}_{r,w}$'s are then obtained through a Laurent expansion of (\ref{Y1k}) in inverse powers of $z$ and $u$. In fact, the expansion  can be derived from the expansion of $z(u-1)/(zu-z-u)$ followed by a shifting of the $w$ index of $u$ by $k$ due to the overall $u^{-k}$ factor in (\ref{Y1k}). Hence,
\begin{align} \label{expand1}
 \frac{z(u-1)}{zu-z-u} & = \left(1-\frac{1}{u} \right)\frac{1}{1-\frac{1}{z}-\frac{1}{u}} = \left(1-\frac{1}{u} \right) \sum_{r=0}^{\infty} \sum_{w=0} ^{\infty } \frac{(r+w)!}{r!w!}\frac{1}{z^r}\frac{1}{u^w} \\ \nonumber
&= \sum_{r=0}^{\infty} \sum_{w=0} ^{\infty } \frac{(r+w)!}{r!w!}\frac{1}{z^r}\frac{1}{u^w} - \sum_{r=0}^{\infty} \sum_{w=1} ^{\infty } \frac{(r+w-1)!}{r!(w-1)!}\frac{1}{z^r}\frac{1}{u^w}\\ \nonumber
&= \sum_{r=0}^{\infty} \sum_{w=0} ^{\infty } \left(\frac{(r+w)!}{r!w!} - \frac{(r+w-1)!}{r!(w-1)!}\right) \frac{1}{z^r}\frac{1}{u^w}  \left(\frac{\mbox{(i) assuming } (-1)!=\Gamma(0)=\infty,} {\mbox{\hspace{-0.8in}(ii) considering (\ref{pa}).} }\right)\\ \nonumber
& =  \sum_{r=0}^{\infty} \sum_{w=0} ^{\infty }  \frac{r(r+w-1)!}{r!w!}  \frac{1}{z^r}\frac{1}{u^w}
\end{align}
where in the third line we have extended the summation from $w=0$ because the denominator of the second term in the bracket vanishes assuming the factorial is analytically continued to a Gamma function such that $n!=\Gamma(n+1)$. Also the term $w=r=0$ in the sum of the last equality should be understood, in the view of (\ref{pb}) and (\ref{pc}),  as taking the limit $w\rightarrow 0$ first followed by $r\rightarrow 0$ yielding a unit coefficient. Generally, the previous arguments employ the ordering of the limits defined in (\ref{pa}) (see also the discussion below (\ref{pa2e})).

Combining now (\ref{expand1}) with (\ref{Y1k}) we obtain
\begin{align}\label{Y1kfinal}
Y^{(k)}(z,u) =  \sum_{r=0}^{\infty} \sum_{w=k} ^{\infty }  \frac{r(r+w-k-1)!}{r!(w-k)!}  \frac{1}{z^r}\frac{1}{u^w}= \sum_{r=0}^{\infty} \sum_{w=0} ^{\infty }  \frac{r(r+w-k-1)!}{r!(w-k)!}  \frac{1}{z^r}\frac{1}{u^w}
\end{align}
where the in the second equality we have again extended the sum from $w=0$ in the view of the fact that $\Gamma(-n)=\infty,$ when $n\,=\,0,\,1,\,2,\,...\,.$
From (\ref{Y1kfinal}) and (\ref{Yfdef}) we identify the coefficients $f^{(k)}_{r,w}$. Using then (\ref{ans1}) we obtain the desired probability formula
\begin{align} \label{p1old}
p^{(k)}(r,w) = \frac{r (r+w-k-1)!w!}{(r+w)!(w-k)!},\, \forall\, r,\,w\geq 0,\, k>0.
\end{align}

\vspace{0.2cm}
\underline{ Case $k= 0$}
\vspace{0.2cm}

In this case, the non-zero terms in the right hand side of (\ref{Y1}) are the $f^{(k=0)}_{r,0} \,\forall \,r  \geq 0$ including the $f^{(k=0)}_{0,0}$ term. In fact, all of these terms are equal to one in the view of (\ref{pbb}); they represent the probability to remain with $k=0$ white balls when we start without any white balls at all and this should happen with probability equal to one. The fact that $f^{(k=0)}_{r,0}=\delta_0^0=1$ has already taken into account the rescaling pre-factor of (\ref{ans1}), which also reduces to one when $w=0$. 

Resumming then the right hand side of (\ref{Y1}) and simplifying yields

\begin{align}\label{Y10}
Y^{(k=0)}(z,u) = \frac{z(u-1)}{zu-z-u}  \mbox{ (because $f^{(k=0)}_{r,0}=1\, \,\forall \,r $ and $f^{(k=0)}_{0,w}=\delta^0_w$).}
\end{align}
Comparing the right hand side of (\ref{Y10}) with that of (\ref{Y1k}) we realize that the probabilities for $k=0$ can be obtained from (\ref{p1old}) by taking the limit $k \to 0$.

To conclude, the general solution is given by

\begin{align} \label{p1}
\boxed{p^{(k)}(r,w) = \frac{r (r+w-k-1)!w!}{(r+w)!(w-k)!}},\, \forall k,\,r,\,w \geq 0,\, k\leq w,
\end{align}
and it satisfies the normalization condition
\begin{align} \label{normp2}
\sum_{k=0}^wp^{(k)}(r,w) =1.
\end{align}
One may check that formula (\ref{p1}) satisfies (\ref{rec1}) and has the right behavior in several limiting cases. Indeed, $p^{(k)}(r,w<k)=0$ where we use the fact that the factorial of negative integer has a first order pole (tends to infinity). Also, in the view of (\ref{pa}), the probability formula yields $p^{(k)}(0,w)=\delta_w^k$.

\subsection{A combinatorial derivation}

In section \ref{2.1} we followed a rather complicated route in order to illustrate part of the formalism of Problem 1. However, Problem 2 may also be solved using simple combinatorics, which we present below for comparison.
 
 We begin by defining the events $A$ and $B$ where
\begin{subequations}
 \begin{align}
 A\,\,= &\mbox{ the event to draw a red ball given that currently $k$ white and one red balls}\\ \nonumber
& \mbox{ remained,}\\
B\,\,= &\mbox{ the event to draw $r-1$ red and $w-k$ white balls out of a population }\\ \nonumber
  &\mbox{ of $r$ red and $w$ white balls.}
 \end{align}
 \end{subequations}
The following two key facts are emphasized: 
(i) The complement of the event $A$,  $A^{\intercal}$, and the event to remain with $k$ white balls starting with $r$ reds and $w$ whites, are mutually exclusive. 
(ii) The event $B$ does not impose any constrains on the possible orderings where the $r-1$ red and the $w-k$ white balls could be chosen.

 With these two key facts at hand we apply the law of total probability by conditioning on the events $A$ and $A^{\intercal}$. Defining $P(A)$ and $P(B)$ as the probabilities for the occurrence of the events $A$ and $B$ respectively we obtain
 \begin{align}\label{p1comb}
 p^{(k)}(r,w)& = p^{(k)}(r,w|A^{\intercal})P(A^{\intercal}) + p^{(k)}(r,w|A)P(A) = 0 + p^{(k)}(r,w|A)P(A)  \\ \nonumber
 & = p(B) \frac{1}{1+k}  =  \frac{\binom{r}{r-1} \binom{w}{w-k}}{\binom{r+w}{r+w-k-1}} \frac{1}{1+k} =  \frac{r (r+w-k-1)!w!}{(r+w)!(w-k)!},
 \end{align}
 which is identical to (\ref{p1}) and this completes the derivation.
 
It is interesting to note that the probability of the event $B$, according to the second equality of the second line of (\ref{p1comb}), essentially follows a type $_2F_1$ hypergeometric distribution\footnote{We remind the reader that ``In probability theory and statistics, the hypergeometric distribution is a discrete probability distribution that describes the probability of $k$ successes (random draws for which the object drawn has a specified feature) in $n$ draws, without replacement, from a finite population of size $N$ that contains exactly $K$ objects with that feature, wherein each draw is either a success or a failure. If a random variable $X$ follows the hypergeometric distribution then its probability mass function is given by $p_X(k) = \frac{\binom{K}{k} \binom{N-K}{n-k} }{\binom{N}{n} }$" (Wikipedia).}. In problem $1$, as we will see further below in sections \ref{genPDE}, \ref{result} and \ref{result4}, the hypergeometric function is a key object.

For later comparison with the results of Problem 1, we provide the probability generating functional of Problem 2 for fixed $r$ and $w$, which is given by 
 \begin{align}\label{Gp2}
 G(r,w;z) \equiv \sum_{w=0}^kp^{(k)}(r,w)z^k = \frac{r}{r+w}  \,_2F_1\left(1,-w,1-r-w ;z\right).
 \end{align}
 It is not a coincidence that equation (\ref{Gp2}) involves a type $_2F_1$ hypergeometric function.

\section{Setting up the problem: From the law of total probability to the required PDE}\label{pde}
We now return to Problem 1. Working as in previous section, we will derive a PDE rather than an algebraic equation for the analogous $Y^{{(k)}}$ we saw earlier.

Using the law of total probability by conditioning on the first removal, we obtain
\begin{align}\label{rec3}
p^{(k)}(r,w) = \frac{r^2}{(r+w)^2}p^{(k)}(r-1,w) + \frac{w^2+2 r w}{(r+w)^2}p^{(k)}(r,w-1)
\end{align}
where the weighting probability of $p^{(k)}(r-1,w)$ is complement to the one of $p^{(k)}(r,w-1)$, and can be understood as $w/(w+r)+ rw/(r+w)^2$: the probability to remove a white ball equals the probability to either chose a white ball directly or to choose a red followed by a white ball.

\subsection{The PDE of the generating functional with rescaling}\label{pdeRes}

In this case too, we make an educated ansatz on (\ref{rec3}) as follows
\begin{align}\label{ans2}
p^{(k)}(r,w) = \left( \frac{r!w!}{(r+w)!}\right)^2 f^{(k)}(r,w)
\end{align}
obtaining the much simpler equation\footnote{Note added after the completion of this work: According to (\ref{psum}), the ansatz $p^{(k)}(r,w) = \frac{(r!)^2 k! w!}{((r+w)!)^2} f^{(k)}(r,w)$ yields an equation of the form (\ref{grec}), which is an even simpler equation for $ f^{(k)}(r,w)$ than equation (\ref{rec4}).}
 for $f^{(k)}$
\begin{align}\label{rec4}
wf^{(k)}(r,w) =  wf^{(k)}(r-1,w) + (2r+w)f^{(k)}(r,w-1).
\end{align}
This ansatz, as in Problem 2, is motivated because the expression $ (r!w!)^2/ ((r+w)!)^2$ provides the probability to remain with $k=w$ white balls; that is if we keep removing only red balls.
Working then as in Problem 1 by multiplying both sides of (\ref{rec4}) by $z^{-r}$ and $u^{-w}$ and summing over $r$ and $w$ from one to infinity, after some algebra, we eventually obtain
\begin{align}\label{sum2}
\frac{1}{z}\sum_{r=0}^{\infty} \sum_{w=0} ^{\infty }wf^{(k)}_{r,w} z^{-r} u^{-w} 
&+\frac{1}{u}\sum_{r=0}^{\infty} \sum_{w=0} ^{\infty }(2r+w+1)f^{(k)}_{r,w} z^{-r} u^{-w} -\sum_{r=0}^{\infty} \sum_{w=0} ^{\infty }wf^{(k)}_{r,w} z^{-r} u^{-w} \\ \nonumber
&= \frac{1}{u} \sum_{w=0} ^{\infty }(w+1)f^{(k)}_{0,w}  u^{-w} -\sum_{w=0} ^{\infty }wf^{(k)}_{0,w}  u^{-w}.  
\end{align}
Defining $Y^{(k)}(z,u)$ as in (\ref{Yfdef}) and making the substitutions
\begin{align}\label{partials}
r\rightarrow -z\partial_z,     \,\,\,\,\, w\rightarrow -u\partial_u
\end{align}
on the $r$ and the $w$ factors that multiply the quantity $f^{(k)}_{r,w} z^{-r} u^{-w}$, equation (\ref{sum2}) yields

\begin{equation}\label{pdeG}
\boxed{
\begin{aligned}
-\left(\frac{u}{z}-u+1\right)\partial_uY^{(k)}(z,u)&-2\frac{z}{u}\partial_zY^{(k)}(z,u)+\frac{1}{u}Y^{(k)}(z,u) \\ 
&=  \frac{1}{u} \sum_{w=0} ^{\infty }(w+1)f^{(k)}_{0,w}  u^{-w} -\sum_{w=0} ^{\infty }wf^{(k)}_{0,w}  u^{-w}.
\end{aligned}
}
\end{equation}
In what follows, we proceed using standard techniques of PDEs. In particular, we proceed using the fact that a first order PDE in two dimensions is equivalent to a system of two first order ODEs. In order to see that we define and solve the first ODE 
\begin{align}\label{dudz}
\frac{du}{dz} = \frac{\frac{u}{z}-u+1}{2\frac{z}{u}} \implies u = \frac{z}{1+z+x_2\sqrt{z}}
\end{align}
where $x_2$ is the integration constant. The next step is to solve for $x_2=x_2(z,u)$ and make the change of variables
 \begin{subequations}\label{x12}
\begin{align}
\label{x1}
x_1&= z,\\
\label{x2}
x_2 &= \frac{z - u - z u}{\sqrt{z} u} = -\sqrt{z}\left(1-\frac{1}{u} + \frac{1}{z} \right)
\end{align} 
\end{subequations}  
where we note that as $z,\,u\to \infty$ we have that $x_1\to \infty$ and $x_2 \to -\sqrt{z} \to -\infty$.  
In the final step, we apply the transformations (\ref{x12}) and the differential equation (\ref{pdeG}) becomes
\begin{equation}\label{odeG}
\boxed{
\hspace{-0.18cm}\begin{aligned}
\frac{1+x_1+ x_2\sqrt{x_1}}{x_1} &\left (Y^{(k)}(x_1,x_2) -2x_1\partial_{x_1}Y^{(k)}(x_1,x_2)  \right)=\\
 \frac{1}{u} \sum_{w=0} ^{\infty }(w+1)f^{(k)}_{0,w}  \left( u(x_1,x_2)\right) ^{-w} &-\sum_{w=0} ^{\infty }wf^{(k)}_{0,w}   \left( u(x_1,x_2)\right)^{-w}, \,\,\,\, u=\frac{x_1}{1+x_1+x_2\sqrt{x_1}},\,\, z=x_1
\end{aligned}
}
\end{equation}
where, as expected, the PDE reduces to a first order ODE. This is the second equation in the equivalent set of the two first order ODEs mentioned earlier.

In the following sections, we will study (\ref{odeG}) in three different versions depending on the right hand side. The first version, section \ref{special}, corresponds to known results in the literature (Problem 3), the second version, sections \ref{genPDE} and \ref{result}, corresponds to the original problem, Problem 1, Rule III, while the third version corresponds to section \ref{result4}, Problem 1, Rule IV.

\subsection{The PDE of the generating functional without rescaling}\label{pdeNoRes}

For completeness, we present the PDE that corresponds to the initial recursion (\ref{rec3}) rather to  (\ref{rec4}), and, which provides the actual probability generating functional. Working as before and defining
\begin{align}\label{Yt}
\tilde{Y}^{(k)}(z,u) = \sum_{r,w=0} ^{\infty } p^{(k)}(r,w)\frac{1}{z^r}\frac{1}{u^w},
\end{align}
the PDE satisfied by $\tilde{Y}^{(k)}(z,u)$ is given by
\begin{align}\label{pde21}
\left(\frac{1}{z}+\frac{1}{u} \right) \tilde{Y}^{(k)} (z,u)
&-\frac{1}{u} \left(u+z(2+u) \right) \partial_z \tilde{Y}^{(k)}(z,u)
-\left(1+u \right) \partial_u\tilde{Y}^{(k)} (z,u)\\ \nonumber
&-z\left(z-1 \right) \partial^2_{zz}\tilde{Y}^{(k)}(z,u)
-2z\left(u-1 \right) \partial^2_{zu}\tilde{Y}^{(k)}(z,u)
-u\left(u-1 \right) \partial^2_{uu}\tilde{Y}^{(k)} (z,u)\\ \nonumber
&=\frac{1}{z}\sum_{r=0}^{\infty} (r+1)^2p^{(k)}(r,0)\frac{1}{z^r}
+\frac{1}{u}\sum_{w=0}^{\infty} (w+1)^2p^{(k)}(0,w)\frac{1}{u^w}\\ \nonumber
&-\sum_{r=0}^{\infty} r^2p^{(k)}(r,0)\frac{1}{z^r}
-\sum_{w=0}^{\infty} w^2p^{(k)}(0,w)\frac{1}{u^w}.
\end{align}

We now have two cases to deal with. We start with Rule III for which conditions (\ref{pb}) and (\ref{pbb}) yield $p_{III}^{(k)}(0,w) =\delta_w^k$ and $p_{III}^{(k)}(r,0) =\delta_0^k$ respectively. Given this and using the fact that
\begin{align}
\frac{1}{z}\sum_{r=0}^{\infty} (r+1)^2p_{III}^{(k)}(r,0)\frac{1}{z^r} -\sum_{r=0}^{\infty} r^2p_{III}^{(k)}(r,0)\frac{1}{z^r} 
= \left(\frac{1}{z}\sum_{r=0}^{\infty} (r+1)^2\frac{1}{z^r} -\sum_{r=0}^{\infty} r^2\frac{1}{z^r} \right) \delta_0^k=0,
\end{align}
we conclude that for the Rule III case, equation (\ref{pde21}) reduces to
\begin{subequations}\label{pde22}
\begin{empheq}[box=\widefbox]{align}
\label{pde3}
 D_{z,u} \tilde{Y}_{III}^{(k)} (z,u)& = (k+1)^2\frac{1}{u^{1+k}} - k^2 \frac{1}{u^k},  \mbox{ for Rule III,}   \\
\label{difO}
D_{z,u}  &\equiv \left(\frac{1}{z}+\frac{1}{u} \right) 
-\frac{1}{u} \left(u+z(2+u) \right) \partial_z 
-\left(1+u \right) \partial_u \\  \notag
&-z\left(z-1 \right) \partial^2_{zz}
-2z\left(u-1 \right) \partial^2_{zu}
-u\left(u-1 \right) \partial^2_{uu}.
\end{empheq}
\end{subequations}
The second case is about Rule IV. Using an analogous approach as in the case of Rule III, it is eventually found that in the case of Rule IV, equation (\ref{pde21}) reduces to

\begin{equation}\label{pde224}
\boxed{
\begin{aligned}
D_{z,u} \tilde{Y}_{IV}^{(k)} (z,u)& = (k+1)^2\frac{1}{z^{1+k}} - k^2 \frac{1}{z^k} ,  \mbox{ for Rule IV}
\end{aligned}
}
\end{equation}
where $D_{z,u}$ is the (same) differential operator defined in (\ref{difO}).


We observe that the resulting PDEs are $2^{nd}$ order linear PDEs with mixed derivatives, which are generally harder equations to solve.

\section{Reproducing known results from the literature}\label{special}
The following known problem \cite{ MuL} is given.

\vspace{0.2cm}
\hspace{-0.6cm}
{\bf Problem 3}
\vspace{0.2cm}

{\it Inside a box containing $r$ red and $w$ white balls random removals occur. The balls are removed one by one according to the three following rules. Rule I: If a white ball is chosen it is immediately discarded from the box. If a red ball is chosen, it is placed back into the box and a new ball is randomly chosen. The second ball is then removed irrespective of the color. Rule II: Once a ball is removed, the game continues from Rule I. Rule III: The game ends once all balls, except from one ball that remains, are removed. The question posed is the determination of the probability that the remaining ball is white.}

\vspace{0.2cm}
\hspace{-0.6cm}
{\bf Solution}
 
\vspace{0.2cm}

All the steps that lead to (\ref{odeG}) are identical for this problem as well and hence, (\ref{odeG}) is the starting point. Given Rule III, which in particular, it implies that even at the event where all reds are removed and where $k>1$ white balls remain, we are still allowed to keep removing white balls (with probability one) until we end up with a single white ball. Thus, we now drop the superscript $k=1$ from $Y^{(k)}$ and $f^{(k)}$ as, according to Rule III, we will keep removing balls until one remains. Therefore, $f_{0,0}=0$ and $f_{0,w}=1, \,\,\forall \,w>0$ \footnote{In particular, equation (\ref{pb}) ceases to apply.} and hence, the right hand side of (\ref{odeG}) becomes $-1/u=-(1+x_1+x_2\sqrt{x_1})/x_1$ (see (\ref{dudz}) and (\ref{x12})). Hence, (\ref{odeG}) simplifies to

\begin{subequations} \label{odespecialTemp}
\begin{empheq}{align}
  Y(x_1,x_2) -2x_1\partial_{x_1}Y(x_1,x_2) & = -1, \\ 
  \lim_{x_2 \to -\frac{1}{\sqrt{x_1}}-\sqrt{x_1}} Y(x_1,x_2) & = 0
\end{empheq}
\end{subequations}
where the boundary condition implies that when no white balls exist ($u\to\infty $), which is equivalent to $x_2 \to -\frac{1}{\sqrt{x_1}}-\sqrt{x_1}$ (see (\ref{x2})), the probabilities and hence, the generating functional to remain with one white ball at the end, should vanish.

The general solution to (\ref{odespecialTemp}) is given by
\begin{subequations} \label{odespecial}
\begin{align}
  Y(x_1,x_2)& = -1+C(x_2)\sqrt{x_1}, \\ 
  \lim_{x_2 \to -\frac{1}{\sqrt{x_1}}-\sqrt{x_1}} Y(x_1,x_2) & = 0
\end{align}
\end{subequations}
where $C(x_2)$ is an arbitrary integration constant to be specified by the boundary condition, which in particular, it implies
\begin{align}\label{C2special}
C(x_2)& = \frac{1}{\sqrt {\frac{1}{2} \left( -2+x_2^2-x_2\sqrt{4-x_2^2} \right)  }} \mbox{, (using $x_2<0)$}\\ \nonumber
&=\frac{1}{\sqrt {\frac{1}{2} \left( -2+x_2^2+x_2^2\sqrt{1-\frac{4}{x_2^2}} \right)  }} \mbox{, (simplifying)} \\ \nonumber
&=\frac{x_2}{2} \left(\sqrt{1-\frac{4}{x_2^2}} -1\right),
\end{align}
where the simplification from the second to the third line uses the following 
\begin{align}\label{idt}
\mbox{if } \sqrt{x}+\frac{1}{\sqrt{x}}=-t \mbox{ then } x = \frac{1}{2} \left( -2+t^2 \pm t^2\sqrt{1-\frac{4}{t^2}}  \right) = \left( \frac{t}{2}\left(1\pm \sqrt{1-\frac{4}{t^2}}  \right)\right)^2.
\end{align}
The last simplification is very convenient when expanding in Laurent series. 

Combining now (\ref{odespecial}) with (\ref{C2special}) we arrive at the desired generating functional
\begin{align}\label{gfs}
\boxed{ Y(x_1(z,u),x_2(z,u)) = -1+\frac{x_2(z,u)}{2} \left(\sqrt{1-\frac{4}{x_2^2(z,u)}} -1\right)\sqrt{x_1(z,u)} }
\end{align}
where $x_1$ and $x_2$, according to (\ref{x12}), are functions of $z$ and $u$. It is straightforward to check that (\ref{gfs}) satisfies the set of equations in (\ref{odespecialTemp}).

The final step is to use (\ref{gfs}) and  (\ref{x12}) and Laurent expand in inverse powers of $z$ and of $u$ noting from (\ref{x12}) that as $z,\,u \to \infty$ we have $x_1 \to \infty$ and $\,x_2 \to -\infty$. In order to expand we make use of the following identities
\begin{subequations}\label{exp}
\begin{align}
\label{expa}
\frac{x_2}{2} \left(\sqrt{1-\frac{4}{x_2^2}} -1\right)&= - \sum_{i=0}^{\infty} \frac{1}{1+i} \frac{\Gamma(1+2i)}{\Gamma^2(1+i)} \left( \frac{1}{x_2}\right)^{1+2i},\\
\label{expb}
\left(\frac{1}{x_2}\right)^l &=  \left(-\sqrt{z}\left(1-\frac{1}{u} + \frac{1}{z} \right)\right)^{-l} \\ \nonumber
&= (-1)^l z^{-\frac{l}{2}} \sum_{n=0}^{\infty}  \sum_{m=0}^{\infty}  \frac{(-1)^m\Gamma(1-l)}{\Gamma(1-l-n-m)\Gamma(1+n)\Gamma(1+m)} \left(\frac{1}{z}\right)^{n} \left(\frac{1}{u}\right)^{m}
 \end{align} 
\end{subequations}
where we use Gamma functions instead of factorials for reasons that will become evident in what follows. Using the expansions (\ref{exp}), the generating functional (\ref{gfs}) expands as
\begin{align}\label{Ysexp1}
 Y(z,u) = -1 +\sum_{i=0}^{\infty}  \sum_{n=0}^{\infty}  \sum_{w=0}^{\infty} &\Bigg[ \frac{1}{1+i} \frac{\Gamma(1+2i)}{\Gamma^2(1+i)}  \times \\ \nonumber 
 &\frac{(-1)^{w+2i+2}\Gamma(-2i)}{\Gamma(-2i-n-w)\Gamma(1+n)\Gamma(1+w)} \left(\frac{1}{z}\right)^{n+i} \left(\frac{1}{u}\right)^{w}\Bigg]\\ \nonumber
 = \sum_{\substack{i,n,w= 0\\(i,n,m)\neq (0,0,0)}}^{\infty}  &\Bigg[\frac{1}{1+i} \frac{\Gamma(1+2i)}{\Gamma^2(1+i)}  \times \\ \nonumber 
 &\frac{(-1)^{w+2i+2}\Gamma(-2i)}{\Gamma(-2i-n-w)\Gamma(1+n)\Gamma(1+w)} \left(\frac{1}{z}\right)^{n+i} \left(\frac{1}{u}\right)^{w} \Bigg].
\end{align}
Changing the dummy indices $i$ and $n$ by setting $i = r-n$ with $r=0,\,1,\,2,\,...$ and (a new) $n=0,\,1,\,2,\,...,\,r$ equation (\ref{Ysexp1}) becomes

\begin{align} \label{fb1}
 Y(z,u) = \sum_{\substack{r,w= 0\\(r,w)\neq (0,0)}}^{\infty}  
 \left(\frac{1}{z}\right)^{r} \left(\frac{1}{u}\right)^{w} \sum_{n=0}^{r} &\frac{1}{1+r-n} \frac{\Gamma(1+2r-2n)}{\Gamma^2(1+r-n)}
\frac{1}{\Gamma(1+n)\Gamma(1+w)}  \times \\ \nonumber
&(-1)^{w}\frac{\Gamma(2n-2r)}{\Gamma(-2r+n-w)}.
\end{align} 
Our approach is to use (\ref{fb1}) and compute the summation term over $n$, which, according to (\ref{Yfdef}), is by definition equal to $f_{r,w}$ (with the $k$ index dropped).  With $f_{r,w}$ at hand, using (\ref{ans2}) we will eventually obtain the required probability formula.

A comment is in order: both, the numerator and the denominator of the last fraction of the Gamma functions in (\ref{fb1}) diverge because $n$, $w$ and $r$ are positive integers with $r \geq n$. This is why we have replaced the factorials with Gamma functions as an analytic continuation of the factorials. These infinities are fictitious and are a consequence of interchanging the summation orderings. Similar infinities will appear in sections \ref{result} and \ref{result4}. Hence, in order to proceed we need the behavior of the Gamma functions near their poles; the corresponding expansions are documented in appendix \ref{A}. As we will show, the end result, when taking the limits carefully, will be finite and will reproduce exactly the correct coefficients such that the recursion (\ref{rec3}) and the boundary condition $p(r,0)=0$ are both fulfilled, which is what matters at the end. 

Using (\ref{ga}) in order to cancel the poles of the Gamma functions\footnote{We only consider the poles resulting from $r\to \mathbb{Z}_{\geq 0}$ because we ignore poles resulting from the dummy summation index $n\to \mathbb{Z}_{\geq 0}$. Another way to see it is to only consider the pole of the variable $-(r-n) \leq 0$. Furthermore, we ignore the pole in $\Gamma(-2(r-n)-n-w)$ coming from $w\to \mathbb{Z}_{\geq 0}$ (and from $n\to \mathbb{Z}_{\geq 0}$) because it is not compensated by a similar $w$ (or $n$) pole in the numerator and hence, it does not contribute a finite part.} the last fraction of the Gamma's in (\ref{fb1}) is exchanged according to
\begin{align}
(-1)^{w}\frac{\Gamma(2n-2r)}{\Gamma(-2r+n-w)}=(-1)^{w}\frac{\Gamma(-2(r-n))}{\Gamma(-2(r-n)-n-w)} \to (-1)^{n}\frac{\Gamma(2r+w-n+1)}{\Gamma(1+2r-2n)}. 
\end{align}
Simplifying then the $\Gamma(1+2r-2n)$ that appears in the numerator and in the denominator of the resulting step, keeping track of the correct overall sign as $(-1)^n$, and changing the dummy index $n$ setting $n\to r-n$, transforms the summand of (\ref{fb1}) into $(-1)^{r+n} \Gamma(1+r+w+n)/\left((1+n)\Gamma(1+w) \Gamma(1+r-n) \Gamma^2(1+n) \right)$ where as before the (new) $n=0,\,1,\,...,\,r$.
Performing then the summation of the resulting expression over $n$ using the summation identity
\begin{align} 
\sum_{n=0}^r (-1)^{r+n}\frac{ \Gamma(1+r+w+n)} { (1+n)\Gamma(1+w) \Gamma(1+r-n) \Gamma^2(1+n) } =\left(\frac{(r+w)!}{r! w!}\right)^2 \frac{w}{(r+1)(r+w)}
\end{align}
proved in Appendix \ref{C}, equation (\ref{fb1}) eventually becomes\footnote{In Appendix \ref{B}, we will briefly provide a second derivation of (\ref{ffinal2}).}
\begin{align}\label{ffinal2}
Y(z,u) = \sum_{\substack{r,w= 0\\(r,w)\neq (0,0)}}^{\infty}  
 \left(\frac{1}{z}\right)^{r} \left(\frac{1}{u}\right)^{w}\left(\frac{(r+w)!}{r! w!}\right)^2 \frac{w}{(r+1)(r+w)}.
\end{align} 
The last step is to use (\ref{Yfdef}) in order to identify the $f_{r,w}$ coefficient from (\ref{ffinal2}) and multiply it by the overall ansatz factor from (\ref{ans2}) to obtain
\begin{subequations}\label{solp2}
\begin{empheq}[box=\widefbox]{align}
  p(0,0)& = 0,  \\
p(r,w)& = \frac{w}{(w+r)(1+r)},\,\, \forall\, r,\,w \,\geq0,\,(r,w)\neq(0,0).
\end{empheq}
\end{subequations}
Equation (\ref{solp2}) is the main result of the section and reproduces the result of \cite{ MuL} obtained through a completely different (inductive) approach. 
\section{Deriving the generic generating functional for any number of remaining white balls $k$, Rule III}\label{genPDE}

Having solved Problem 3, which uses all the ingredients of our proposed approach, we gained confidence in order to attack Problem 1. In the following two sections, we simplify the notation by dropping the index $_{III}$ from $Y_{III}^{(k)}(x_1,x_2)$, $f^{(k)}_{III;r,w}$ and $p_{III}^{(k)}(r,w)$; we will restore it at the very end on the final formula for $p_{III}^{(k)}(r,w)$.

Starting from (\ref{odeG}), we write the differential equation and the boundary conditions for Problem 1, Rule III. Noting that $f^{(k)}_{0,w} =\delta^k_w$, because $f^{(k)}_{0,w}$ denotes the probability to remain with $k$ white balls if we start with no red balls and with $w$ white balls (see (\ref{pb})), we obtain
\begin{subequations} \label{odep1}
\begin{empheq}[box=\widefbox]{align}
\label{odep1a}
  \frac{1+x_1+ x_2\sqrt{x_1}}{x_1} \Big(Y^{(k)}(x_1,x_2)& -2x_1\partial_{x_1}Y^{(k)}(x_1,x_2)  \Big) \\ \nonumber
  &= (1+k)\frac{1}{u^{1+k}(x_1,x_2)}-k\frac{1}{u^{k}(x_1,x_2)},\\
\label{odep1b}
\lim_{z\to = \infty}Y^{(k)}(x_1,x_2) &= \frac{1}{u^k},\,\,z=z(x_1, x_2),\, u=u(x_1,x_2)
\end{empheq}
\end{subequations}
where we have indicated explicitly that $z$ and $u$  are functions of $x_1$ and of $x_2$ (see (\ref{x12})).

The boundary condition $z\to \infty$ (game starts without any red balls) is a boundary case where the game ends before it (even) begins. Then, the probability to remain with $k$ white balls should be equal to one if we start with $k$ white balls and zero otherwise \footnote{The argument takes into account the ansatz factor in (\ref{ans2}), which also reduces to one, when $r=0$.}. The boundary condition in the generating functional representation, equation (\ref{odep1b}), is a manifestation of equation (\ref{pb}), which is the boundary condition in the probability representation.

The solution of (\ref{odep1}) that satisfies the boundary condition turns out to be
\begin{subequations} \label{solk}
\begin{empheq}[box=\widefbox]{align}
\label{solka}
 & \hspace{-0.0in}Y^{(0)}(x_1,x_2)  = 1 -\sqrt{x_1}  \frac{2x_2+\sqrt{2\left(x_2^2 \left(1+\sqrt{1-\frac{4}{x_2^2}}\,\right)-2\right)}}{x_2^2\left(1+\sqrt{1-\frac{4}{x_2^2}} \,\right)-4} \\ \nonumber
& \hspace{0.8in}=1 + \sqrt{x_1}  \frac{x_2\left(\sqrt{1-\frac{4}{x_2^2}} -1\,\right)}{x_2^2\left(1+\sqrt{1-\frac{4}{x_2^2}} \,\right)-4}, \\
 \label{solkb} &\hspace{-0.2in}Y^{(k>0)}(x_1,x_2)  = \frac{2^{-2(1+k)}x_1^{-k}}{x_2^2\left(1-\frac{4}{x_2^2} \right)} \left(2+x_2 \sqrt{x_1}\left(1-\sqrt{1-\frac{4}{x_2^2}}\,\right) \right)^{1+2k} \times \\ \nonumber
& \hspace{-0.35in}  \Bigg \{ \hspace{-0.1in} - \hspace{-0.05in} \frac{2+x_2\sqrt{x_1}}{1+x_1+x_2\sqrt{x_1}} \hspace{-0.06in}   \left(\hspace{-0.05in} 2\hspace{-0.03in} +\hspace{-0.03in} x_2 \sqrt{x_1}\left(1\hspace{-0.05in} + \hspace{-0.05in} \sqrt{1-\frac{4}{x_2^2}}\,\right) \hspace{-0.05in} \right) 
\hspace{-0.05in}   \Bigg(\hspace{-0.05in} 1\hspace{-0.03in} + \hspace{-0.03in} \frac{2x_2\sqrt{x_1}\sqrt{1-\frac{4}{x_2^2}}}{2+x_2 \sqrt{x_1}\left(1-\sqrt{1-\frac{4}{x_2^2}}\,\right)} \Bigg)^k \\ \nonumber
&\hspace{-0.3in}   \hspace{0.25in}+ 2\frac{(1+k)x_2^2-2 }{1+2k} \,_2F_1\left(-1-2k,-k;-2k;-\frac{2x_2\sqrt{x_1}\sqrt{1-\frac{4}{x_2^2}}}{2+x_2 \sqrt{x_1}\left(1-\sqrt{1-\frac{4}{x_2^2}}\,\right)} \right) \Bigg \}
\end{empheq}
\end{subequations}
where in (\ref{solka}) we simplified using (\ref{idt}), and where $_2F_1$ in (\ref{solkb}) is a hypergeometric function.

It is a straightforward but also a tedious process to show that (\ref{solk}) satisfies the differential equation (\ref{odep1}) and a computer program such as ``Mathematica" is recommended. The fact that $Y^{(0)}$ fulfills the boundary condition can be verified by taking first the limit $\lim {x_2 \to -\sqrt{z}}(1-\frac{1}{u})$ followed by  the limit $\lim {z \to \infty}$. These limits are a consequence of (\ref{x12}) and of the boundary condition (\ref{odep1b}).

In order to see that the boundary condition for $Y^{(k>0)}$ is also satisfied we take the same limits as the limits we took for $Y^{(0)}$. Then, the overall multiplicative factor of (\ref{solkb}) decays as $-(u-1)^{-3-2k}u^2/2\times1/z^{1+k}$. Expanding the first term of the curly bracket (not the one involving the $_2F_1$) to leading order in $z$ shows that this quantity grows as $ -2(u-1)^{2k+2}/u^{k+1} \times z^{k+1}$. Next, we move to the second term in the curly bracket, that is the hypergeometric function whose argument grows as $-(1-u)^2/u \times z$. Taking into account that if $k$ is a positive integer, we find that the $_2F_1\left(-1-2k,-k;-2k;-x\right)$ becomes a terminating polynomial of degree $k$, which grows as  $(2k+1)/(k+1)x^{k}$ as $x \to \infty$. Thus, the whole term involving the $_2F_1$ grows as $2(1-u)^{2k+2}/u^{k+2}  \times z^{k+1}$. Therefore, the curly bracket grows as $-2(u-1)^{2k+3}/u^{k+2} \times z^{k+1}$ and given the overall pre-factor $-(u-1)^{-3-2k}u^2/2\times1/z^{1+k}$ it is concluded that (\ref{odep1b}) converges to $1/u^k$ as should. We have thus just shown that the boundary condition is fulfilled $\forall \,k>0$. 

In the following section we will finalize the solution of Problem 1, Rule III. Using (\ref{solk}), we will compute explicitly the probabilities $p^{(k)}(r,w)$ for the first few values of $k$ and from there we will derive the formula for any $k$.

It is also interesting to note that the alternative boundary conditions
\begin{align}\label{x2bc}
\lim_{x_2\to -\sqrt{x_1}-\frac{1}{\sqrt{x_1}}}Y^{(k)}(x_1,x_2) &= \delta_0^k \frac{x_1}{x_1-1}
\end{align}
yield the same exact solutions provided by equations (\ref{odep1}). In particular, the boundary conditions defined by (\ref{x2bc}) are the continuum version of the discrete boundary conditions provided by equation (\ref{pbb}). 

These boundary conditions imply that as $u\to \infty$ (game starts without any white balls), which due to (\ref{x12}) is equivalent to $x_2\to -\sqrt{x_1}-\frac{1}{\sqrt{x_1}}$, the probability to remain with $k>0$ white balls is zero. If, on the other hand, $k=0$ then the probability to remain with no white balls starting without any white balls, must be equal to one independently on the initial number of red balls \footnote{The argument takes into account the ansatz factor in (\ref{ans2}), which also reduces to one, when $w=0$.}. Indeed, expanding the expression $x_1/(x_1-1)$ in inverse powers of $x_1$ yields unit coefficients \footnote{We remind the reader that $x_1=z$ and that the variable $z$ corresponds to the red balls.}. We reserve the details, which show that the boundary conditions (\ref{x2bc}) are satisfied by (\ref{odep1}), for appendix \ref{Z}.

 In case the reader wonders about the uniqueness of the problem, the answer is that the two boundary conditions are equivalent because they define the same exact problem.  In particular, one may choose to evolve the initial data starting from the surface $z \to \infty$ or to evolve the initial data starting from the surface $u \to \infty$, as long as the chosen boundary surface is known. In this problem, it happens that we know both boundaries. In what follows, things will become clearer and the simultaneous fulfillment of both boundary conditions will be checked once explicit formulas are obtained.

\section{The general probability formula for Rule III}\label{result}

Having developed all the necessary machinery we are now in the position to answer the question of Problem 1, Rule III. The strategy we will apply is as follows. We first consider (\ref{solk}) and take the limits $k \to 0,\,1\,,\,2$. Then, noting from (\ref{x12}) that as $z,\,u \to \infty$ we have $x_1\to \infty$ and $x_2\to -\infty$,  we Laurent expand the solutions at infinity, that is in inverse powers of $x_1=z$ and $x_2$. Then, we re-expand $x_2$ in inverse powers of $z$ and $u$ in an analogous way as we worked in Problem 3, section \ref{special}. The final step is to identify the coefficients of $1/z^r$ and $1/u^w$ as the $f_{r,w}^{(k)}$'s of (\ref{Yfdef}) and from there, using (\ref{ans2}), we will obtain the $p^{(k)}(r,w)$. Using then the explicit expressions for $p^{(0)}(r,w)$, $p^{(1)}(r,w)$ and $p^{(2)}(r,w)$ and by studying the pattern, we will make an ansatz for the generic solution $p^{(k)}(r,w)$ for any $k$. The last step is to check that the ansatz fulfills the initial recursion we begun with, namely equation (\ref{rec3}), and the condition  $p^{(k)}(0,w) = \delta_w^k$. It is also noted that in principle we could had reached the same answer without making any anstaz at all. That is through an even more systematic approach; by Laurent expanding the general PDE solutions (\ref{solk}) for arbitrary $k$.

\subsection{The case $k=0$}\label{sk0}
Despite this case can be obtained using the strategy outlined in the introduction of the current section, there is a way around it using the results of section \ref{special}. In particular, the probability of Problem 1, Rule III for $k=0$ is equal to the complement of the probability of Problem 3, which is the probability that the last ball is red; that it no white balls remain at the end of the game. Thus, using (\ref{solp2}), we find
\begin{subequations}\label{solp10temp}
\begin{empheq}{align}
  p^{(0)}(0,0)& = 1,  \\
p^{(0)}(r,w)& = r\frac{1 + r + w}{(w+r)(1+r)},\,\, (r,w)\neq(0,0).
\end{empheq}
\end{subequations}
In the view of (\ref{pa}) and (\ref{pc}), equation (\ref{solp10temp}) can be written compactly as
\begin{align}\label{solp10}
\boxed{
p_{III}^{(0)}(r,w) = r\frac{1 + r + w}{(w+r)(1+r)},\,\, \forall \,r,\,w \geq 0.
}
\end{align}
As a cross check of (\ref{solp10}), we substitute (\ref{x12}) into (\ref{solka}) and expand the result at $z=\infty$ and at $u=\infty$ up to $4^{th}$ order to obtain
\begin{align}
Y^{(0)}(z,u) = \left(1+O\left( \frac{1}{u^5}\right) \right) &+ \left(1+ \frac{3}{u}+\frac{6}{u^2}+ \frac{10}{u^3}+\frac{15}{u^4}+O\left( \frac{1}{u^5}\right) \right) \frac{1}{z}\\ \nonumber\
&+ \left(1+ \frac{8}{u}+\frac{30}{u^2}+ \frac{80}{u^3}+\frac{175}{u^4}+O\left( \frac{1}{u^5}\right) \right) \frac{1}{z^2}\\ \nonumber\
&+ \left(1+ \frac{15}{u}+\frac{90}{u^2}+ \frac{350}{u^3}+\frac{1050}{u^4}+O\left( \frac{1}{u^5}\right) \right) \frac{1}{z^3}\\ \nonumber\
&+ \left(1+ \frac{24}{u}+\frac{210}{u^2}+ \frac{1120}{u^3}+\frac{4410}{u^4}+O\left( \frac{1}{u^5}\right) \right) \frac{1}{z^4}+O\left( \frac{1}{z^5}\right).
\end{align}
We then compare previous expansion with the quantity
\begin{align}
Y^{(0)}(z,u) = \sum_{r=0}^{4} \sum_{w=0} ^{4 } f^{(0)}_{r,w} z^{-r} u^{-w}  = 1+\mathop{\sum_{r=0}^{4}\sum_{w=0}^{4}}_{(r,w)\neq(0,0)} \left(\frac{(r+w)!}{r! w!}\right)^2 r\frac{1 + r + w}{(w+r)(1+r)} z^{-r} u^{-w} 
\end{align}
observing an exact matching order to order in $z$ and in $u$. It is also notable that both boundary conditions (\ref{odep1b}) and (\ref{x2bc}) are satisfied. Indeed, checking the condition (\ref{odep1b}) is obvious. 
On the other hand, as $u\to \infty$, the expansion reduces to $1+1/z+1/z^2+1/z^3+....=z/(z-1)=x_1/(x_1-1)$ in accordance with the (alternative) boundary condition in (\ref{x2bc}). Thus, as discussed below (\ref{x2bc}), we can see explicitly the equivalence of the two boundary conditions.

Equation (\ref{solp10}) will also be obtained as the $k\to 0$ limit of the general solution $p^{(k)}(r,w)$. 

\subsection{The case $k=1$} \label{sk1}

Taking the limit in (\ref{solkb}) as $k \to 1$ yields
\begin{align}\label{k1}
Y^{(1)}(x_1,x_2) &=\left(1 +\frac{2}{3x_1} + \frac{x_2}{\sqrt{x_1}}  \right) \\ \nonumber
&+\frac{x_2\sqrt{x_1}}{6}\left(3-\sqrt{1-\frac{4}{x_2^2}}\,\, \right)-\frac{x_2^3\sqrt{x_1}}{6}\left(1-\sqrt{1-\frac{4}{x_2^2}}\,\, \right).
\end{align}
Using (\ref{x12}), the first line of (\ref{k1}) expands as
\begin{align}\label{3z}
1 +\frac{2}{3x_1} + \frac{x_2}{\sqrt{x_1}}  = \frac{1}{u}-\frac{1}{3 z}
\end{align}
while the second line, after some algebra, expressed as an expansion in $x_1$ and $x_2$ yields
\begin{align}\label{intk1}
\frac{x_2\sqrt{x_1}}{6}\left(3-\sqrt{1-\frac{4}{x_2^2}}\,\, \right)&-\frac{x_2^3\sqrt{x_1}}{6}\left(1-\sqrt{1-\frac{4}{x_2^2}}\,\, \right) \\ \nonumber
& = -\sqrt{x_1}\sum_{i=2}^{\infty} \frac{(i-1)\Gamma(2i-1)}{i(1+i)\Gamma^2(i)} \frac{1}{x_2^{2i-1}}.
\end{align}
From here and on we work analogously to equations (\ref{exp})-(\ref{fb1}). In particular\footnote{The explanations that follow assume an identical indexing in $i$, $n$ $r$ and $w$ as the indexing of section \ref{special}.}, we expand $ \frac{1}{x_2^{2i-1}}$ using (\ref{expb}), exchange the index $i$ with $i=r-n+1$ where $r=1,\,2,\,...$ and $n=0,\,1,\,2,\,...\,r-1$  as before to finally obtain 
\begin{align}\label{prevk11}
Y^{(1)}(z,u)= \sum_{r=1}^{\infty} \sum_{w=0}^{\infty}  \left(\frac{1}{z}\right)^{r} \left(\frac{1}{u}\right)^{w} 
\Bigg \{ \sum_{n=0}^{r-1} &\frac{(r-n)}{(1+r-n)(2+r-n)} \frac{\Gamma(1+2(r-n))}{\Gamma^2(1+r-n)}
\frac{(-1)^{w}}{\Gamma(1+n)\Gamma(1+w)}   \\ \nonumber
&\times \frac{\Gamma(2(n-r))}{\Gamma(-2r+n-w)} \Bigg \}+ \frac{1}{u}-\frac{1}{3z}.
\end{align} 
We note that for $r\geq1$, $n<r$ and $w \geq0$ the numerator and the denominator in the last fraction of the Gammas in (\ref{prevk11}) diverges. Using then (\ref{gb}), which allows us to replace the last fraction of the Gamma functions with $(-1)^{n+w}\Gamma(1+2r-n+w)/\Gamma(1+2(r-n))$ \footnote{We only consider the poles resulting from $r\to \mathbb{Z}_{\geq 0}$ because we ignore poles resulting from the dummy summation index $n\to \mathbb{Z}_{\geq 0}$. Furthermore, we ignore the pole in $\Gamma(-2r+n-w)$ coming from $w\to \mathbb{Z}_{\geq 0}$ because it is not compensated by a similar $w$ pole in the numerator and hence, it does not contribute a finite part.}, and equation (\ref{Yfdef}), which defines the $f^{(k)}_{r,w}$ coefficients, we obtain
\begin{align}
\label{prevk12}
f^{(1)}_{r,w} &=
\sum_{\substack{n= 0\\r \geq1}}^{r-1}  
\frac{(-1)^n(r-n)}{(1+r-n)(2+r-n)} \frac{1}{\Gamma^2(1+r-n)}
\frac{\Gamma(1+2r-n+w)}{\Gamma(1+n)\Gamma(1+w)} +\left(\frac{1}{u}-\frac{1}{3z} \right) \Bigg|_{r,w} \\ \nonumber\
&=\sum_{\substack{n= 1\\r \geq1}}^{r}  \frac{(-1)^{r-n} n}{(1+n)(2+n)} \frac{1}{\Gamma^2(1+n)}
\frac{\Gamma(1 + r + w + n)}{\Gamma(1+r-n)\Gamma(1+w)}  +\left(\frac{1}{u}-\frac{1}{3z} \right) \Bigg|_{r,w}, 
\,\,r\geq 1, 
\,\,w\geq 0
\end{align} 
where in the second equality we changed the dummy index according to $n\to r-n$ and where $ |_{r,w}$ denotes the coefficient of $\frac{1}{z^r}\frac{1}{u^w}$ in the preceding bracket.
Equation (\ref{prevk12}) can be partitioned in four cases according to
\begin{subequations} \label{prevk13}
\begin{align}
\label{prevk13a}
f^{(1)}_{1,0} & = \frac{1}{3}+\left( \frac{1}{u}-\frac{1}{3z} \right) \Bigg|_{1,0}=0, \\
\label{prevk13d}
f^{(1)}_{0,1} & =\left( \frac{1}{u}-\frac{1}{3z} \right) \Bigg|_{0,1}=1,\\
\label{prevk13b}
f^{(1)}_{r,w} &= \frac{(r+w-2)!(1+r+w)!}{(r-1)!(2+r)!(w-1)! w!}, \mbox{ $r\geq1,\,w\geq0$ and $(r,w)\neq(1,0)$}\\
f^{(1)}_{r,w} &=0 \mbox{ otherwise}
\end{align}
\end{subequations}
where (\ref{prevk13a}) is obtained from (\ref{prevk12}) for $(r,w)=(1,0)$, and it has exactly the form we need in order to cancel the undesired $1/(3z)$ term coming from (\ref{prevk11}). Equation (\ref{prevk13b}) is obtained by performing the summation using the identity
\begin{align}\label{bio2}
\frac{(r+2)!}{(r+w+1)!}\sum_{n=1}^{r}\frac{(-1)^n n}{(1+n)(2+n)} \frac{(r+w+n)!}{(r-n)! (n!)^2} &= \sum_{n=1}^{r} (-1)^n \binom{r+2}{n+2} \binom{r+w+n}{n-1}\\ \nonumber\
&= (-1)^r \frac{(r+w-2)!}{(r-1)!(w-1)!},
\end{align}
which can be proved working along the lines of equations (\ref{C1})-(\ref{CendProof}) with minor modifications. We note that in the limit $(r,w)\to(0,1)$, the right hand side of equation (\ref{prevk13b}) tends to one\footnote{The limit $r=0,\,w=1$ in (\ref{prevk13b}) must be taken according to equation (\ref{pa2e}) (see also comments that follow (\ref{pa2e})). In particular, when taking the first limit $w \to 1$, both of the $(r-1)!$ terms that appear in the numerator and the denominator of (\ref{prevk13b}) cancel out yielding $f^{(1)}_{0,1}=1$.}, and hence, (in theory) it provides the $1/u$ term (see boundary conditions, (\ref{odep1b})); in reality, the $1/u$ term is actually coming from $\left(\frac{1}{u}-\frac{1}{3z} \right) \big|_{r=0,w=1}$. Thus, equation (\ref{prevk13d}) can be absorbed in (\ref{prevk13b}) in the view of the limits ordering of (\ref{pa}), and by extending the range of validity suitably. 

Collecting all the terms of (\ref{k1}) using (\ref{prevk13}) and the observations we just made above, we finally obtain
\begin{align}\label{k1exp}
Y^{(1)}(z,u) = \sum_{r,w=0}^{\infty} f^{(k)}_{r,w}\left(\frac{1}{z}\right)^{r}  \left(\frac{1}{u}\right)^{w}=
\sum_{r=0}^{\infty}  \sum_{w=0}^{\infty} \frac{(r+w-2)!(1+r+w)!}{(r-1)!(2+r)!(w-1)! w!}  \left(\frac{1}{z}\right)^{r} \left(\frac{1}{u}\right)^{w} 
\end{align}
where we extend the summations on $r$ and on $w$ according to $r,\,w \geq0$ in the view of the fact that $1/(r-1)!$ and $1/(w-1)!$ tend separately to zero as $r\to0$ and $w\to0$ respectively. As noted earlier, when $(r,w)=(0,1)$ we get $f^{(k)}_{0,1}=\delta_1^k$ (see (\ref{pa2e}) and the discussion below  (\ref{pa2e}) in order to see how these limits should be taken). The reason we expect $f^{(1)}_{0,1}=1$ is because this coefficient encodes the $1/u$ term (see (\ref{odep1b}) and (\ref{3z})) and hence, equation (\ref{k1exp}) provides the right coefficient with value equal to one, which is precisely the probability to remain with one white ball if we start with one white ball and no red balls. To rephrase, if we start with no red balls $(r=0)$ and $w \neq1$ white balls the probability to remain with $k=1$ white balls is zero, otherwise if $w=1$ the probability is equal to one, and this is precisely the meaning of $f^{(k)}_{0,1}=\delta_1^k$.

The last step is to use (\ref{Yfdef}) in order to identify $f^{(k)}_{r,w}$ from (\ref{k1exp}), multiply by the coefficient of (\ref{ans2}) and do the necessary simplifications to obtain

\begin{align}\label{solp11}
\boxed{
 p_{III}^{(1)}(r,w) = \frac{rw(1 + r + w)}{(1+r)(2+r)(w+r-1)(w+r)},\,\,\forall \, r,\,w \geq0,
 }
\end{align}
which also includes the cases $w=0,\,\, r\neq1$ and $(r,w)=(0,1)$ in the view of (\ref{pbb}) and (\ref{pb}) respectively.
One may verify that equation (\ref{solp11}) satisfies the initial recursion we begun with, namely equation (\ref{rec3}) and the boundary condition (\ref{pb}), and in fact, all the equations (\ref{pb})-(\ref{pe}).

Comparing (\ref{solp10}) with (\ref{solp11}) we already start to see a pattern forming. In the next section, once the $k=2$ case is computed, the pattern will become obvious.

\subsection{The case $k=2$}\label{sk2}

Taking the limit in (\ref{solkb}) as $k \to 2$ yields
\begin{align}\label{k2}
Y^{(2)}(x_1,x_2) &=\left(1 +\frac{3}{5x_1^2} +\frac{3x_2}{2x_1^{\frac{3}{2}}} + \frac{4+3x_2^2}{3x_1}+ \frac{2x_2}{\sqrt{x_1}}  \right) \\ \nonumber
&+\sqrt{x_1}x_2 \left(\frac{1}{2}-\frac{2}{15}\sqrt{1-\frac{4}{x_2^2}}\,\, \right)
+x_2^3\sqrt{x_1}\left(-\frac{1}{3}+\frac{7}{30}\sqrt{1-\frac{4}{x_2^2}}\,\, \right)\\ \nonumber
&+\frac{\sqrt{x_1}x_2^5}{20}\left( 1-\sqrt{1-\frac{4}{x_2^2}}\,\, \right).
\end{align}
Before computing the general $f^{(k)}_{r,w}$ term, we expand (\ref{k2}) up to $4^{th}$ order in $1/z$ and $1/u$ as a way to cross-check our calculations along the way. The result is
\begin{align}\label{k2exp}
Y^{(2)}(z,u)& = \left(\frac{1}{u^2}+O\left( \frac{1}{u^5}\right) \right) + \left(\frac{1}{u^2}+ \frac{5}{3u^3}+\frac{5}{2u^4}+O\left( \frac{1}{u^5}\right) \right) \frac{1}{z}\\ \nonumber\
&+ \left(\frac{1}{u^2}+ \frac{4}{u^3}+\frac{21}{2u^4}+O\left( \frac{1}{u^5}\right) \right) \frac{1}{z^2}
+ \left(\frac{1}{u^2}+ \frac{7}{u^3}+\frac{28}{u^4}+O\left( \frac{1}{u^5}\right) \right) \frac{1}{z^3}\\ \nonumber\
&+ \left(\frac{1}{u^2}+ \frac{32}{3u^3}+\frac{60}{u^4}+O\left( \frac{1}{u^5}\right) \right) \frac{1}{z^4}+O\left( \frac{1}{z^5}\right).
\end{align}
As expected, according to the boundary condition (\ref{odep1b}), the lowest order in $1/u$ is $1/u^2$, which basically says that unless we start with at least two white balls, the probability to end up with two white balls is zero. More specifically, we know that if $(r,w)=(0,2)$ we should get a probability equal to one, which is precisely equal to the coefficient of the $1/z^01/u^2$ term in the expansion. 
We also note that as $u\to \infty$, $Y^{(2)}(z,u) \to 0$ in agreement with the alternative boundary condition (\ref{x2bc}).

Returning now to equation (\ref{k2}) and expanding it yields
\begin{align}\label{y21}
Y^{(2)}(z,u) = \frac{1}{u^2}-\frac{1}{6z}-\frac{1}{2uz}+\frac{1}{10z^2} -2\sqrt{x_1(z,u)}\sum_{i=0}^{\infty}\frac{4^i i \Gamma(\frac{1}{2}+i)}{\sqrt{\pi}\Gamma(4+i)}\frac{1}{x_2^{1+2i}(z,u)}
\end{align}
where the sum corresponds to the last two lines of (\ref{k2}). We note that we already have the desired $1/u^2$ term.
Working then analogously to sections \ref{special} and \ref{sk1}, and using the Gamma function identities
\begin{subequations} \label{idsk2}
\begin{align}
\Gamma \left(\frac{1}{2}+x\right) = 2^{1-2x}\sqrt{\pi}  \frac{\Gamma(2x)}{\Gamma(x)},\\
\pi \csc \left(2\pi x\right) = - 2x \Gamma(-2x) \Gamma(2x),
\end{align}
\end{subequations}
and expanding $x_2$ as in (\ref{expb}) we find that (\ref{y21}) yields
\begin{align} \label{sumk21}
Y^{(2)}(z,u)  &= \frac{1}{u^2}-\frac{1}{6z}-\frac{1}{2uz}+\frac{1}{10z^2} 
 -\sum_{r,w=0}^{\infty}  \left(\frac{1}{z}\right)^{r}  \left(\frac{1}{u}\right)^{w}  \Bigg \{2 (-1)^w  \\ \nonumber
&\times \sum_{n=0}^r \frac{ \pi \csc(2\pi(r-n))}{ \Gamma(1+w) \Gamma(r-n) \Gamma(4+r-n) \Gamma(1+n) \Gamma(-2r-w+n)} \Bigg \}.
\end{align}
We note that for any $w$, when $r=0$, and hence $n=0$, the summation term over $n$ becomes trivial and also it yields a zero result. The reason is because the numerator has only a single pole in $r$ coming from $\csc(2\pi (r-n)) $ while the denominator has a second order pole in $r$, one coming from $\Gamma(r-n)$ and one coming from $\Gamma(-2r-w+n)$. Hence, the summation over $r$ should begin from $r=1$. 

As it turns out the summation over $n$ is, up to overall Gamma function factors, a regularized hypergeometric function with integer coefficients that reduces to a hypergeometric function $_2F_1$. In particular, 
(\ref{sumk21}) yields
\begin{align} \label{sumk22}
Y^{(2)}(z,u)  = & \frac{1}{u^2}-\frac{1}{6z}-\frac{1}{2uz}+\frac{1}{10z^2}  \\ \nonumber
&-\sum_{r=1,w=0}^{\infty} \left(\frac{1}{z}\right)^{r}  \left(\frac{1}{u}\right)^{w}  \Bigg \{ 2 (-1)^w  \frac{ \pi \csc(2\pi r) \,_2F_1\left( -3-r,1-r;-2r-w;1\right)}{ \Gamma(1+w) \Gamma(4+r) \Gamma(r) \Gamma(-2r-w)} \Bigg \}
\\ \nonumber  
=& \frac{1}{u^2}-\frac{1}{6z}-\frac{1}{2zu}+\frac{1}{10z^2} \\ \nonumber
&+\sum_{r=1,w=0}^{\infty} \left(\frac{1}{z}\right)^{r}  \left(\frac{1}{u}\right)^{w} \Bigg \{2   \frac{\Gamma(1+2r+w)  \,_2F_1\left( -3-r,1-r;-2r-w;1\right)}{ \Gamma(1+w) \Gamma(4+r) \Gamma(r) } \Bigg \}
\end{align}
where in the second equality we used equation (\ref{gb}) on $\Gamma(-2r-w)$ and the fact that as $r\to j \in \mathbb{Z}$ then $\pi \csc(2\pi r) \approx \frac{1}{2(r-j)}$. These two facts allow the cancellation of the poles in the ratio $\pi \csc(2\pi r)/\Gamma(-2r-w)$ and its replacement with $-\Gamma(2r+w+1)(-1)^w$ \footnote{We only consider the pole resulting from $r\to \mathbb{Z}_{\geq 0}$ and we ignore the pole in $\Gamma(-2r-w)$ coming from $w\to \mathbb{Z}_{\geq 0}$ because this pole is not compensated by a similar $w$ pole in the numerator and hence, it does not contribute a finite part.}. Given also that $_2F_1\left( -4,0;-2;x\right)=1$, $_2F_1\left( -4,0;-3;x\right)=1$, $_2F_1\left( -5,-1;-4;x\right)=1-5/4 x$ and $_2F_1\left( -5,-1;-5;x\right)=1- x$, the curly bracket in (\ref{sumk22}) yields the following terms
\begin{subequations} \label{cancellterms}
\begin{align}
\frac{1}{6z} &\mbox{ for $(r,w)=(1,0)$,}\\
\frac{1}{2zu} &\mbox{ for $(r,w)=(1,1)$,}\\
-\frac{1}{10z^2} &\mbox{ for $(r,w)=(2,0)$,}\\
0 &\mbox{ for $(r,w)=(2,1)$,}\\
0 &\mbox{ for $r=0$ and $w \geq 0$},\\
\label{w2r2}
0 &\mbox{ for $w<2$ and $r > 2$}.
\end{align}
\end{subequations}
Equation (\ref{w2r2}) is not obvious but it will become soon due to equation (\ref{casesk2b}), which implies that for $w=\{0,1 \}$ and $r>2$, the $\Gamma(w-1)$ term in the denominator diverges without any compensating factor coming from the numerator; thus the whole term is equal to zero. It is noted that the first three sub-equations of (\ref{cancellterms}) have precisely the right form in order to cancel the unnecessary terms from (\ref{sumk22}) (see (\ref{k2exp})), which thus, in the view of (\ref{Yfdef}), yields
\begin{subequations} 
\begin{align} 
f^{(2)}_{r,w}  \left(\frac{1}{z}\right)^{r} \left(\frac{1}{u}\right)^{w} & = \frac{1}{u^2}, \,\, (r,w)=(0,2),\\
f^{(2)}_{r,w}  \left(\frac{1}{z}\right)^{r} \left(\frac{1}{u}\right)^{w} &=2 \left(\frac{1}{z}\right)^{r} \left(\frac{1}{u}\right)^{w}  \frac{\Gamma(1+2r+w) }{ \Gamma(1+w) \Gamma(4+r) \Gamma(r) } \\ \nonumber
&\times \,_2F_1\left( -3-r,1-r;-2r-w;1\right), \,\,\, r \geq 1,\,\, w\geq 2,\\
f^{(2)}_{r,w}  \left(\frac{1}{z}\right)^{r} \left(\frac{1}{u}\right)^{w} &= 0 \mbox{ otherwise}.
\end{align}
\end{subequations}
Next, we use the Gauss theorem on the hypergeometric function $_2F_1$, which states that 
\begin{align}
_2F_1(a,b;c;1) = \frac{\Gamma(c) \Gamma(c-a-b)}{ \Gamma(c-a) \Gamma(c-b)},\,\,  c-a-b>0.
\end{align}
Strictly speaking in our case $c-a-b = 2 -w\leq 0$ for $w\geq2$ and hence, the Gamma functions, of both, the numerator and the denominator, will diverge. Hence, we take the limits of the Gamma functions carefully by taking care of the poles that appear using (\ref{gpoles})\footnote{While taking the required limits, we cross-check by matching the $f^{(2)}_{r,w}$'s with the coefficients of (\ref{k2exp}).}. The end result is finite and is given by
\begin{subequations} \label{casesk2}
\begin{align} 
f^{(2)}_{r,w}  & = 1, \,\, (r,w)=(0,2),\\
 \label{casesk2b}
f^{(2)}_{r,w} &  = 2 \frac{\Gamma(r+w-2) \Gamma(r+w+2)}{ \Gamma(w-1) \Gamma(w+1) \Gamma(4+r) \Gamma(r) }, \,\,\, r \geq 1,\,\, w\geq 2.\\
f^{(2)}_{r,w}  &= 0 \mbox{ otherwise}.
\end{align}
\end{subequations}
Equation (\ref{casesk2}) reproduces the $\frac{1}{z^r} \frac{1}{u^w}$ coefficients of (\ref{k2exp}), and this provides confidence that the equation is correct; it can be written compactly as
\begin{align} \label{f2comp}
f^{(2)}_{r,w} &  = 2 \frac{\Gamma(r+w-2) \Gamma(r+w+2)}{ \Gamma(w-1) \Gamma(w+1) \Gamma(4+r) \Gamma(r) }, \,\, \forall\, r,\,w \geq0
\end{align}
where we have extended the range of $r$ and of $w$ by observing that: (i) If $r=0$ and $w=2$, taking the limits in (\ref{f2comp}) according to (\ref{pa}), yields $f^{(2)}_{0,2}=1$ in agreement with (\ref{pb}). (ii) Else if $r=0$ and $w>2$ we have one $r$ pole in the denominator coming from $\Gamma(r)$ that is not compensated by a similar pole in the numerator and hence, the result is zero. 
(iii) Else if $r=0$ and $w=0,\,1$, despite the ratio $\Gamma(r+w-2)/\Gamma(r)$ is finite, the $1/\Gamma(w-1)$ term tends to zero and hence, $f^{(2)}_{0,1}=f^{(2)}_{0,2}=0$. (iv) Else, if $w=0,\,1$ and $r \geq 1$ 
the result is also zero, including the cases $(r,w)=\{(1,0),\,(1,1)\}$, in the view of the orderings of the limits of (\ref{pa}) and the fact that the $1/\Gamma(w-1)$ in (\ref{f2comp}) goes to zero as $w\to0,\,1$. Given that we are dealing with $k=2$ and $w=0,\,1$, we see that case (iv) is consistent with equation (\ref{pe}) as should.


The last step is to multiply (\ref{f2comp}) by the coefficient of (\ref{ans2}) and do the necessary simplifications to obtain

\begin{align}\label{solp12}
\boxed{
 p_{III}^{(2)}(r,w) =2 \frac{rw(w-1)(1 + r + w)}{(1+r)(2+r)(r+3)(w+r-2)(w+r-1)(w+r)},\,\, r,w \geq 0.
 }
\end{align}
It can be checked that equation (\ref{solp12}) satisfies the initial recursion (\ref{rec3}) and the boundary condition (\ref{pb}), and in fact, all the equations (\ref{pb})-(\ref{pe}).
This completes the derivation of $ p_{III}^{(2)}(r,w)$. In the following section, we will provide the general solution $ p_{III}^{(k)}(r,w)$ $\forall k$.


\subsection{The general solution}

Having computed $p^{(k)}(r,w)$ for $k=0,\,1$ and $2$, we are now in position to attempt for a general ansatz solution. Observing equations (\ref{solp10}), (\ref{solp11}) and (\ref{solp12}) one may guess that the solution should have the form
\begin{align}\label{finallypka}
p_{III}^{(k)}(r,w) &\sim r \frac{w(w-1)(w-2)...(w-k+1)}{(r+1)(r+2)...(r+k+1)} \frac{1+r+w}{(r+w)(r+w-1)....(r+w-k)}\\ \nonumber\
&\sim r\frac{r! w!}{(r+k+1)!(w-k)!} \frac{(1+r+w)(r+w-k-1)!}{(r+w)!}              
\end{align}
up to an overall constant that does not depend on $r$ and $w$. The constant is then specified by the requirement that $p^{(k)}(0,k)=1$ in the view of the boundary condition (\ref{pb}) from where we infer that the solution must be 
\begin{align}\label{finallypk}
\boxed{
p_{III}^{(k)}(r,w) =\frac{k! r!(r+w+1)}{(r+k+1)!}\,\, \frac{rw!(r+w-k-1)!}{(r+w)!(w-k)!}\,\, \forall r,\,w,\,k\geq0,\, k\leq w.
}        
\end{align}

The final step is to verify that (\ref{finallypk}) is the required probability solution. Indeed, it is a matter of straightforward algebra to show that (\ref{finallypk}) satisfies the initial recursion we begun with, namely equation (\ref{rec3}), and all the equations (\ref{pb})-(\ref{pe}).

We do the following two cross-checks: (i) We first note that the solution reproduces the cases $k=0,\,1,\,2$, equations (\ref{solp10}), (\ref{solp11}) and (\ref{solp12}). (ii) Also, as it is shown in Appendix \ref{D}, the probability formula is normalized and hence, it satisfies
\begin{align}\label{normalp}
\sum_{k=0}^w p_{III}^{(k)}(r,w) =1,
\end{align}   
which basically says that when the game ends, we will surely end up with a number of white balls between zero and the initial number $w$. We note that the normalization equation (\ref{normalp}) comes out automatically and this key fact serves as another cross check of the correctness of (\ref{finallypk}). At this stage, the solution is considered as complete.

Two important observations can be made: (a) It is notable that (\ref{finallypk}) is, up to the rescaling factor $\frac{k! r!(r+w+1)}{(r+k+1)!}$, the same as the solution of Problem 2, equation (\ref{p1}). This implies that this overall factor, in a sense, encodes the additional complication of adding red balls back into the box. (ii) One can also show that the (true) generating probability functional $\tilde{Y}^{(k)}(z,u)= \tilde{Y}^{(k)}_{III}(z,u)$ defined in (\ref{Yt}) satisfies the PDE (\ref{pde22}) with $p^{(k)}(r,w)=p_{III}^{(k)}(r,w)$ given by the probability formula (\ref{finallypk}). 

We also compute the probability generating functional for fixed $r$ and $w$. It is is given by

 \begin{align}\label{Gp1III}
 G_{III}(r,w;z) \equiv \sum_{w=0}^kp_{III}^{(k)}(r,w)z^k = \frac{1+r+w}{1+r}  \frac{r}{r+w}  \,_3F_2\left(1,1,-w;2+r,1-r-w ;z\right),
 \end{align}
which involves a $\,_3F_2$ hypergeometric function rather than a $\,_2F_1$, which was the case for Problem 2 (see (\ref{Gp2})).

\subsubsection{Maximal probability}\label{max3}


In this section we investigate the maxima of $p_{III}^{(k)}(r,w)$ as a function of $k$ for fixed values of $r$ and of $w$. We start by constraining $k$ solving the two inequalities
\begin{subequations} 
\begin{align} 
\label{inq2a}
p_{III}^{(k)}(r,w) \geq p_{III}^{(k+1)}(r,w), \\
\label{inq2b}
p_{III}^{(k)}(r,w) \geq p_{III}^{(k-1)}(r,w),
\end{align}
\end{subequations}
which provide the set of $k$'s for which $p_{III}^{(k)}(r,w)$ has local maxima (if any). The inequalities (\ref{inq2a}) and (\ref{inq2b}) imply 
\begin{subequations} 
\begin{align} 
k& \leq \frac{1}{2} (1+r)(r+w)-1, \\
k & \geq \frac{1}{2} (1+r)(r+w)
\end{align}
\end{subequations}
respectively. Evidently the system of inequalities has no solution. This hinds that, for $r>0$, the maximal probability is at the boundary cases $p_{III}^{(0)}(r>0,w)$ or $p_{III}^{(w)}(r>0,w)$. In fact, we guess that the $p_{III}^{(0)}(r>0,w)$ should be the required maximal probability. In order to show the claim, we take the difference $p_{III}^{(k)}(r,w)-p_{III}^{(1+k)}(r,w)$ and simplify to obtain
\begin{align} 
p_{III}^{(k)}(r,w)& - p^{(1+k)}(r,w) \\ \nonumber
&= |C|r\Gamma(r+w-k-1) (-2 - 2 k + (1+r)(r+w)), \, r \geq 1,\,k\leq w-1
\end{align}
where $|C|$ is a positive constant of ratios of factorials. The difference $p_{III}^{(k)}(r,w) - p_{III}^{(1+k)}(r,w)$ is minimum when $k=w-1$, and given that $r\geq1$, it implies that  the difference is always positive and hence, $p_{III}^{(k)}(r,w)>p_{III}^{(1+k)}(r,w)$. Therefore, when $r>0$ we find that $p_{III}^{(k)}(r,w)$ is a monotonically decreasing function of $k$. On the other hand, we know that for $r=0$, $p_{III}^{(k)}(0,w) =\delta_w^k$.
Thus, we have just proved
\begin{subequations} \label{maxpk}
\begin{empheq}[box=\widefbox]{align}
\label{maxpk3a}
&sup\left(p_{III}^{(k)}(r=0,w)|k=0,\,1,\,...\,w \right)=w,\\
\label{maxpk3b}
&sup\left(p_{III}^{(k)}(r\geq1,w)|k=0,\,1,\,...\,w \right)=0. 
\end{empheq}
\end{subequations}
To conclude, the maximal probability occurs for: (i) $k=w$ if $r=0$ and is given by $p_{III}^{(k=w)}(0,w)=1$. (ii) Else, the maximal probability occurs for $k=0$ and the corresponding probability $p_{III}^{(k=0)}(r,w)$ is given by equation (\ref{solp10}).

\subsubsection{Limiting cases}

It is also interesting to investigate the behavior of (\ref{finallypk}) in the limits $r\to \infty$ for $w$ fixed and any $k$, and for $w \to \infty$ for $r$ and $k$ fixed. For this purpose we use the more convenient equation (\ref{finallypka}) from where it is deduced that
\begin{subequations}  \label{limcases}
\begin{align} 
\label{limcasesr}
p_{III}^{(k)}(r,w) &= k!\frac{1}{r^{2k} } \left( \delta_0^k + \prod_{i=0}^{k-1} (w-i) \right)+ O\left(\frac{1}{r^{2k+1}} \right),\,\, r \to \infty, \, r \gg w,\\
\label{limcasesw}
p_{III}^{(k)}(r,w) &=k! \frac{r \,r!}{(r+k+1)! }+ O\left(\frac{1}{w} \right),\,\, w \to \infty, \, w \gg r,\, k,
\end{align}
\end{subequations}
where $\prod_{i=0}^{k-1} $ for $k=0$ is defined to be zero.

The asymptotic expansion (\ref{limcasesr}) says that in the limit $r \to \infty$ with $w$ kept fixed, the probability $p_{III}^{(k>0)}(r,w)$ decays as $\frac{1}{r^{2k}}$. Moreover, if $k=0$, the probability $p_{III}^{(k=0)}(r,w)$ tends to one, which means that when the red balls are much more than the white balls, the game will (most likely) end without any white balls. Both of these two asymptotic results behave as expected. On the other hand, the asymptotic expansion (\ref{limcasesw}) in the limit $w \to \infty$ with $r$ and $k$ kept fixed yields less expected results. In particular, the probability $p_{III}^{(k)}(r,w)$ (to leading order in $w$) becomes independent on the initial number of the white balls $w$.

\section{The general probability formula for Rule IV}\label{result4}

In this section we solve Problem 1, Rule IV. The ideas are similar as those of Problem 1, Rule III and hence, the derivations are sketchy and much shorter compared to those of sections \ref{genPDE} and \ref{result}. In particular, the proofs between several steps are similar or even identical to the aforementioned sections and therefore, they are omitted. 

\subsection{The PDE, the boundary conditions and the solution, Rule IV}
The starting point is the PDE (\ref{odeG}) where the coefficients in the sums in the right hand side of the equation are all equal to $f_{0,w}^{(k)}=\delta_0^k$ (see (\ref{pbb4})). The (physical) reason is because $f_{0,w}^{(k)}$ is (proportional to) the probability to remain with $k$ red balls if we start without any and this should yield probability equal to one if $k=0$ and equal to zero otherwise. Summing then the two series over $w$ we find that the result is zero $\forall k$, including $k=0$. Hence, the PDE reduces to
\begin{subequations}\label{pdebc4}
\begin{align}\label{pde4}
Y^{(k)}(x_1,x_2) -2x_1\partial_{x_1}Y^{(k)}(x_1,x_2) & =0,\\
\label{bc4}
\lim_{x_2 \to -\sqrt{x_1}-\frac{1}{\sqrt{x_1}}} Y^{(k)}(x_1,x_2) & = \frac{1}{x_1^k}
\end{align}
\end{subequations}
where the boundary condition (\ref{bc4}) is explained as follows. The condition $x_2 \to -\sqrt{x_1}-\frac{1}{\sqrt{x_1}}$, according to (\ref{x12}), is equivalent to the condition that all the white balls are out of the box ($u \to \infty$). Then, in this case, the only solution for $k$ red balls to remain is if inside the box exist exactly $k$ red balls. In other words, the lowest order term in the Laurent expansion of $Y^{(k)}(x_1,x_2) $ should be $\sim 1/z^k$. The solution then to (\ref{pdebc4}) is given by\footnote{Analogously to the Rule III case, we could had imposed the alternative boundary conditions $\lim_{x_1 \to \infty} Y^{(k)}(x_1,x_2) =\frac{u}{u-1}\delta_0^k$ and still obtain the same $Y^{(k)}(x_1,x_2)$ given by equation (\ref{sol4}).}

\begin{align}\label{sol4}
Y^{(k)}(x_1(z,u),x_2(z,u)) =  \sqrt{x_1} \left( \frac{x_2}{2} \left( \sqrt{1-\frac{4}{x_2^2}} - 1\right) \right)^{2k+1}.
\end{align}

\subsection{The probabilities for the first few cases, Rule IV}
In this section we Laurent expand (\ref{sol4}) for $k=0,\,\,1,$ and $2$ and obtain explicit formulas for $p_{IV}^{(k)}(r,w)$. In particular, we first expand in inverse powers of $x_2$ and then re-expand in inverse powers of $z$ and of $u$ identifying the coefficients with the $f^{(k)}_{r,w}$'s (see (\ref{ans2})).

\subsubsection{The case $k=0$}

In fact, the solution to this case is already found in section (\ref{special}) and is (almost) given by equation (\ref{solp2}). The precise answer is 
\begin{align}\label{solp04}
\boxed{ p_{IV}^{(0)}(r,w) = w\frac{1}{1+r} \, \frac{1}{w+r}, \forall \, r,\,w \geq 0}
\end{align}
where we extend the applicability of the formula in order to include $(r,w)=(0,0)$ in the view of (\ref{pc4}) and the discussion below equation (\ref{pa2e4}). One can check that equation (\ref{solp04}) satisfies the recursion (\ref{rec3}), and also fulfills the boundary condition (\ref{pb4}).

\subsubsection{The case $k=1$}

It can be shown that (\ref{sol4}) for $k=1$ expands as
\begin{align}
Y^{(1)}(x_1(z,u),x_2(z,u))  &= \sqrt{x_1} \left( \frac{x_2}{2} \left( \sqrt{1-\frac{4}{x_2^2}} - 1\right) \right)^{3} \\ \notag
&= -3 \sqrt{x_1} \sum_{i=0}^{\infty} \frac{i \,\Gamma(2i+1)}{(i+1)(i+2)\Gamma^2(i+1)} \frac{1}{x_2^{2i+1}}.
\end{align}

The next step is to follow exactly the same steps that led from equation (\ref{intk1}) to equation (\ref{k1exp}) to eventually obtain 

\begin{align}\label{solp14}
\boxed{ p_{IV}^{(1)}(r,w) =3 w \frac{ r }{(r+1)(r+2)}\, \frac{ (r+w+1)}{(r+w-1)(w+r)}, \forall \, r,\,w \geq 0.}
\end{align}

One can check that equation (\ref{solp14}) satisfies the recursion (\ref{rec3}), and also fulfills the boundary condition (\ref{pb4}).

\subsubsection{The case $k=2$}

It can be shown that (\ref{sol4}) for $k=2$ expands as
\begin{align}
Y^{(2)}(x_1(z,u),x_2(z,u))  &= \sqrt{x_1} \left( \frac{x_2}{2} \left( \sqrt{1-\frac{4}{x_2^2}} - 1\right) \right)^{5} \\ \notag
&= -\frac{5}{2} \sqrt{x_1} \sum_{i=0}^{\infty} \frac{(i-1)\, i\, 2^{2i+1}\Gamma \left(i+\frac{1}{2}\right)}{\sqrt{\pi} \Gamma(i+4)}\frac{1}{x_2^{2i+1}}.
\end{align}

The next step is to work as in the $k=1$ case by following exactly the same steps that led from equation (\ref{intk1}) to equation (\ref{k1exp}) to eventually obtain 

\begin{align}\label{solp24}
\boxed{ p_{IV}^{(2)}(r,w) =5 w \frac{(r-1)r }{(r+1)(r+2)(r+3)} \, \frac{(r+w+1)(r+w+2)}{(r+w-2)(r+w-1)(w+r)}, \forall \, r,\,w \geq 0.}
\end{align}
One can check that equation (\ref{solp24}) satisfies the recursion (\ref{rec3}), and also fulfills the boundary condition (\ref{pb4}).

\subsection{The general solution}

Equations (\ref{solp04}), (\ref{solp14}) and (\ref{solp24}) motivate the following ansatz for the general solution for any k
\begin{align}\label{finallypka4}
p_{IV}^{(k)}(r,w) = (2k+1)w \frac{(r-k+1)...(r-1)r }{(r+1)(r+2)...(r+k+1)} \, \frac{(r+w+1)(r+w+2)...(r+w+k)}{(r+w-k)...(r+w-1)(w+r)},
\end{align}
which can be re-written as
\begin{align}\label{finallypk4}
\boxed{p_{IV}^{(k)}(r,w) =(2k+1)\frac{ r!(r+w+k)!}{(r+k+1)!(w+r)!}\,\, \frac{wr!(r+w-k-1)!}{(r+w)!(r-k)!},\,\, \forall\, r,\,w,\,k\geq0,\, k\leq w.}
\end{align}
One may check that (\ref{finallypk4}) satisfies both, the recursive equation (\ref{rec3}) and the boundary condition (\ref{pb4}), and in fact, all the consistency-check equations (\ref{pb4})-(\ref{pe4}). As a cross check, we verified that equation (\ref{finallypk4}) reproduces (\ref{solp04}), (\ref{solp14}) and (\ref{solp24}) and most importantly, it satisfies the normalization condition
\vspace{-0.2in}
\begin{align}\label{normalp4}
\sum_{k=0}^r p_{IV}^{(k)}(r,w) = 1.
\end{align}
We note that the normalization equation (\ref{normalp4}) comes out automatically. This key fact serves as another cross check of (\ref{finallypk4}). At this stage, the solution is considered as complete.

Two important observations can be made: (a) It is notable that (\ref{finallypk4}) is, up to the rescaling factor $\frac{(2k+1) r!(r+w+k)!}{(r+k+1)!(w+r)!}$, the same as the solution of Problem 2, equation (\ref{p1}) with the (expected) reflection $ r \leftrightarrow  w$. This implies that this overall factor, in a sense, encodes the additional complication of adding red balls back into the box. (ii) One can also show that the (true) generating probability functional $\tilde{Y}^{(k)}(z,u) = \tilde{Y}^{(k)}_{IV}(z,u)$ defined in (\ref{Yt}) satisfies the PDE (\ref{pde224}) with $p^{(k)}(r,w)=p_{IV}^{(k)}(r,w)$ given by the probability formula (\ref{finallypk4}). 

We conclude the section by providing the probability generating functional 
\vspace{-0.1in}
 \begin{align}\label{Gp1IV} \hspace{-0.03cm}
 \hspace{-0.1cm}G_{IV}(r,w;\hspace{-0.06cm}z)& \equiv \sum_{r=0}^kp_{IV}^{(k)}(r,w)z^k = \frac{w}{(1+r)(r+w)} \Big[ \,_3F_2\left(1,-r,1+r+w;2+r,1-r-w ;z\right)\\ \notag
 &+ 2r \frac{1+r+w}{(2+r)(r+w-1)} z \,_3F_2\left(2,1-r,2+r+w;3+r,2-r-w ;z\right)\Big]
 \end{align}
for fixed $r$ and $w$. The generating functional involves a linear combination of $\,_3F_2$ hypergeo- metric functions rather than a single $\,_3F_2$, which was the case for Problem 1, Rule III (see (\ref{Gp1III})).

\subsubsection{Maximal probability}
This section investigates the set of $k$'s for which $p_{IV}^{(k)}(r,w)$ is maximized.
For the boundary case $w=0$ we have $p_{IV}^{(k)}(r,0)=\delta_r^k$ (see (\ref{pb4})) and hence, the probability is maximized for $k=r$ with maximal value $p_{IV}^{(r)}(r,0)=1$.
The other boundary case is when $r=0$ in which case $p_{IV}^{(k)}(0,w)=\delta_0^k$ (see (\ref{pbb4})) and hence, the probability is maximized for $k=0$ with maximal value $p_{IV}^{(0)}(0,0)=1$.
Otherwise, if $w>0$ and $r>0$, we work analogously to section \ref{max3} and we find that the set of $k$'s that maximize (locally) the $p_{IV}^{(k)}(r,w)$ is specified by the two inequalities
\begin{subequations} 
\begin{align} 
k&\geq \sqrt{\frac{(1+r)(r+w)}{2w-1}}-1, \,r,\,w \geq1, \\
k & \leq \sqrt{\frac{(1+r)(r+w)}{2w-1}},\,r,\,w \geq1.
\end{align}
\end{subequations} 
We observe that the inequalities are linear in $k$, which means that (modulus degeneracies) there exists only one unique maximum. 
Before investigating the maxima, we define $k_0$ by
\begin{align} \label{k0}
k_0 \equiv \sqrt{\frac{(1+r)(r+w)}{2w-1}}.
\end{align}
Given that the investigation here is for $r \geq1$ and $w \geq1$, it is deduced that $k_0>1$. We now want to constrain $k_0$ from above. As it can be shown, there are two cases: Case 1. If $w=1$, then $k_0=r+1$ and hence, the maximum occurs for $k=r$. Case 2. Else if $w>1$, then $k_0<r+1$. Case 2 is partitioned into two sub-cases: Case 2a. If $k_0 \in \mathbb{Z}_+$ then there are two maxima that correspond to $k=k_0-1$ and to $k=k_0$, which occur with equal probabilities. Case 2b. Otherwise if $k_0 \notin \mathbb{Z}_+$ then there exists a unique maximum that corresponds to the integer between $k_0 -1$ and $k_0$; that is for $k= \left \lfloor{k_0}\right \rfloor $ (the integer part of $k_0$). This is always the case for $r=1$ and $w>1$; the maximum occurs for $k=1$ because $k_0 \notin \mathbb{Z}_+$ with $k_0 \in (1,2)$.

Collecting all cases together, we have just proved

\begin{subequations} \label{maxpk4}
\begin{empheq}[box=\widefbox]{align}
\label{maxpk4a}
&sup\left(p_{IV}^{(k)}(r,w=0)|k=0,\,1,\,...\,r \right)=r,\\
\label{maxpk4b}
&sup\left(p_{IV}^{(k)}(r=0,w)|k=0 \right)=0,\\
\label{maxpk4c}
&sup\left(p_{IV}^{(k)}(r \geq 1,w=1)|k=0,\,1,\,...\,r \right)=r,\\
\label{maxpk4d}
&sup\left(p_{IV}^{(k)}(r\geq1,w \geq2)|k=0,\,1,\,...\,r \right)=\{k_0-1,k_0 \}, k_0\in \mathbb{Z}_+,\\
\label{maxpk4e}
&sup\left(p_{IV}^{(k)}(r \geq1,w \geq2)|k=0,\,1,\,...\,r \right)=  \left \lfloor{k_0}\right \rfloor, \, k_0 \notin \mathbb{Z}_+.
\end{empheq}
\end{subequations}
This ends our investigation of the maxima of $p_{IV}^{(k)}(r ,w )$.

\subsubsection{Limiting cases}

In this section we investigate the behavior of (\ref{finallypk4}) in the limits $w\to \infty$ for $r$ fixed and any $k$, and for $r \to \infty$ for $w$ and $k$ fixed. For this purpose we use the more convenient equation (\ref{finallypka4}) from where it is deduced that
\begin{subequations}  \label{limcases4}
\begin{align} 
\label{limcasesr4}
p_{IV}^{(k)}(r,w) &= (2k+1)  \frac{w}{r^2}+ O\left(\frac{1}{r^{3}} \right)= (2k+1)p_{IV}^{(0)}(r,w) + O\left(\frac{1}{r^{3}} \right),\,\, r \to \infty, \, r \gg w,\,k,\\
\label{limcasesw4}
p_{IV}^{(k)}(r,w) &=(2k+1)\frac{(r!)^2}{(r-k)!(r+k+1)!}+ O\left(\frac{1}{w} \right),\,\, w \to \infty, \, w \gg r.
\end{align}
\end{subequations}

The asymptotic expansion (\ref{limcasesr4}) says that in the limit $r \to \infty$ with $w$ and $k$ kept fixed, the probability $p_{IV}^{(k)}(r,w)$ decays as $\frac{1}{r^{2}}$ independently on $k$ \footnote{Modulus the $2k+1$ overall factor. More precisely, in this limit, the probabilities grow linearly in $k$ according to $(2k+1)p_{IV}^{(0)}(r,w)$; they grow as integer multiples of a fundamental quantity, the probability $p_{IV}^{(0)}(r,w)$.} contrary to the case of Rule III, equation (\ref{limcasesr}). Moreover, the asymptotic expansion (\ref{limcasesw4}) in the limit $w \to \infty$ with $r$ kept fixed is also interesting. In particular, the probability $p_{IV}^{(k)}(r,w)$ (to leading order in $w$) becomes independent on the initial number of the white balls $w$ analogously to the behavior of the Rule III case, equation (\ref{limcasesw}). In particular, if $r=1$, then $p_{IV}^{(0)}(1,w)\approx p_{IV}^{(1)}(1,w) \approx \frac{1}{2}$ for $w \gg 1$.

\section{Discussion}\label{diss}

Wow!! That has been a long but joyful journey for such a rather simple problem. One could think that our approach amounts in shooting a bug with a bazooka! While such a view might be right, in this work we have seen a concrete example of how a recursion equation arising from a probability problem, without a systematic method to approach, was transformed into a differential equation that was solved using standard PDE methods. Hence, we have seen a concrete example of how difference equations, such as those arising from discrete random processes, and for which the tools in the literature are less developed, could be transformed into differential equations, where the literature is rich and well studied. In particular, the approach we followed has been algorithmic and it may be summarized by the following recipe: 

\vspace{0.2in}
\hspace{0.4cm} {\bf Step 1.} Usage of the total probability law in order to derive a difference equation (rec-

\hspace{2.1cm} ursion)\hspace{-0.04cm} in\hspace{-0.04cm} one \hspace{-0.04cm}or \hspace{-0.04cm}more \hspace{-0.04cm}variables.\hspace{-0.04cm} Impose \hspace{-0.04cm}suitable \hspace{-0.04cm}discrete \hspace{-0.04cm}boundary \hspace{-0.04cm}conditions. 

\vspace{0.2in}
\hspace{0.4cm} {\bf Step 2.} Apply z-transformations in order to transform the recursion into a differential 

\hspace{2.1cm} equation for the probability generating functional. 

\vspace{0.2in}
\hspace{0.4cm} {\bf Step 3.} Translate the discrete level boundary conditions of Step 1 into continuum 

\hspace{2.1cm} boundary conditions at the (P)DE level.

\vspace{0.2in}
\hspace{0.4cm} {\bf Step 4.} Solve the (P)DE enforcing the boundary conditions of step 3. 

\vspace{0.2in}
\hspace{0.4cm} {\bf Step 5.} Laurent expand the solutions taking care of any fictitious infinities that may 

\hspace{2.1cm}appear by taking suitable limits. Hence, obtain the required probabilities.

\vspace{0.2in}

This recipe is adaptable to large classes of probability problems, and generally for problems involving difference equations. For instance, adapting the ideas of Problem 1, Rule III, to those of Rule IV, whose results have been obtained fast and straightforwardly, has been effortless.


Retracing the steps (for Problem 1, Rule III), we started from a probability question and derived a recursion, equation (\ref{rec3}), with two variables using the law of total probability. Then, through a 2D $z-$transformation (\ref{Yfdef}), we eventually derived a first order 2D PDE (\ref{pdeG}) with suitable source terms and boundary conditions (\ref{odep1}) whose solution, up to a factor (see (\ref{ans2})), provides the generating probability functional of the problem at hand for any number of remaining white balls $k$. Given that a first order 2D PDE is equivalent to a set of two first order ODEs, the PDE, through suitable coordinate transformations (see (\ref{x12})), was reduced into two decoupled first order ODEs, equations (\ref{dudz}) and (\ref{pdeG}). This step, in a sense, decouples the 2D problem into two 1D problems of suitable scaling variables, which are functions of the initial ones. While the two scaling variables are decoupled at the DE level, they are yet coupled in a non-trivial way through the boundary conditions (i.e. see (\ref{solk})). Solving then the differential equation, and hence, obtaining the general generating functional (\ref{solk}), we performed a closed form Laurent expansion. Using the formula from the expansion, which provides the coefficients $f^{(k)}(r,w)$ (see (\ref{Yfdef}) and (\ref{ans2})), we were eventually able to find the general probability formula, equation (\ref{finallypk}). Finally, we checked that the probability formula satisfies the initial recursion equation we begun with, and its boundary conditions (\ref{pb}), thus completing the solution of the problem. An analogous approach has been followed in Problem 1, Rule IV.

It is notable how a relatively simple probability problem to state and to understand, has  involved such a heavy computational machinery. The most challenging step has been to Laurent expand the generating functional and to deal with the fictitious infinities that appeared on the way due to the interchanging of the summation orderings. These infinities showed up as poles of Gamma functions; this part has been considerably more involved than solving the PDE itself. 

We have also seen that the probability generating functionals involve hypergeometric functions (see (\ref{Gp1III})) or linear combinations of hypergeometric functions (see (\ref{Gp1IV})) of type $_3F_2$. These results generalize the simpler version of the problem, which does not require placing balls back into the box, and, which involves a $_2F_1$ distribution instead (see (\ref{Gp2})). Hence, in a sense, the additional complication of adding balls back into the box is captured by extending the $_2F_1$ type of the probability generating functional into a $_3F_2$ type.

As a bonus, we found the power series solution of the rather complicated family of PDEs given by equation (\ref{pde22}), which satisfy the same boundary conditions $\tilde{Y}_{III}^{(k)}(z \to \infty,u) =1/u^k$ as the boundary conditions satisfied by $Y_{III}^{(k)}(z,u) $ (see (\ref{odep1b})). The solution is given by $\tilde{Y}_{III}^{(k)}(z,u) = \sum_{r,w=0}^{\infty}p_{III}^{(k)}(r,w) \frac{1}{z^r} \frac{1}{u^w} $. Likewise, for Rule IV, $\tilde{Y}_{IV}^{(k)}(z,u) = \sum_{r,w=0}^{\infty}p_{IV}^{(k)}(r,w) \frac{1}{z^r} \frac{1}{u^w} $ is the solution of the PDE (\ref{pde224}), which satisfy the same boundary conditions $\tilde{Y}_{IV}^{(k)}(z,u \to \infty,u) =1/z^k$ as 
the boundary conditions satisfied by $Y_{IV}^{(k)}(z,u) $ (see (\ref{bc4}) and (\ref{x12})).

It would be interesting to find simple probability problems, such as the ones presented in this work, whose differential equation representation is the same as that of known problems from other areas of mathematics or physics or even from finance. That would provide a sort of duality between a physically meaningful and involved problem, and a rather simple probability problem such as a box of balls of various colors.

\section*{Acknowledgments}

I would like to thank A. Kryftis for keep challenging me and in particular, A. Anastasiou (LSE and Un. of Cyprus) for reading the manuscript and for providing insightful comments, and D. Christofides (UCLAN, Cyprus) for providing the solution of appendix \ref{Gminus}, and also S. Agapiou (Un. of Cyprus), S. Hormann (Graz Un. of Technology, Austria), T. Bruss (Universite Libre de 
Bruxelles), and E. Mossel and S. Sheffield (MIT) for referring me to useful sources. I also thank my friends at Vincent House (London), and especially G. Vogiatzi for encouraging me to publish this work, and my colleagues at JPMorgan Chase, and in particular, S. El Hamoui, A. Eriksson, M. Green, J. Lorenzen and S. Mcgarvie for engaging in stimulating discussions\footnote{Disclaimer: \hspace{-0.04cm}The current \hspace{-0.04cm}paper \hspace{-0.04cm}consists of a personal \hspace{-0.04cm}work \hspace{-0.04cm}curried out by \hspace{-0.04cm}A. \hspace{-0.04cm}Taliotis, \hspace{-0.04cm}and it expresses \hspace{-0.02cm}his \hspace{-0.02cm}own personal \hspace{-0.04cm}views on the \hspace{-0.04cm}particular topic. \hspace{-0.04cm}In particular, \hspace{-0.04cm}this  \hspace{-0.04cm}work is not, in  \hspace{-0.04cm}any way, \hspace{-0.04cm}endorsed  \hspace{-0.04cm}by  \hspace{-0.04cm}or  \hspace{-0.04cm}related with \hspace{-0.04cm}JPMorgan \hspace{-0.04cm}Chase (the firm) or its employees or stakeholders or any interests or entities represented by the firm.}. A warm thank you to Staxto for ``being here'', during the long nights after the JPMC office hours, while preparing this write-up. I would like to express my deep gratitude to my teacher and good friend Y. Kovchegov from The Ohio State University for patiently showing me, among other, how to carry on long and tedious calculations and how to, in a quantum field theoretical language, ``renormalize" (deal with) infinities, and eventually obtain finite and meaningful results. Lastly, I would like to thank Nikolas and Leo-Anastasis for invading into my routine teaching me the beauty in life, and for all the time I have taken away from them.


 



\appendix
\renewcommand{\theequation}{A\arabic{equation}}
\setcounter{equation}{0}
\section{Gamma function expansions near poles}
\label{A}

For Gamma function arguments $ \in  0\cup\mathbb{Z}_-$, the following key equation provides the residues 
\begin{align}\label{gres}
Res   \left( \Gamma,-n \in  0\cup\mathbb{Z}_+ \right) = \frac{(-1)^n}{n!}=\frac{(-1)^n}{\Gamma(1+n)}.
\end{align}
Using (\ref{gres}), the following expansions, which are needed in the intermediate \hspace{-0.02cm}steps, can be derived.
\vspace{-0.1in}
\begin{subequations}\label{gpoles}
\begin{align}
\label{ga}
\lim_{x\to n}\Gamma(-x-m)& \approx \frac{(-1)^{1+ n+m}}{\Gamma(1+n+m) } \frac{1}{x-n} +O(1),\mbox{$\,\,m,\,n \in 0\cup \mathbb{Z}^+$,}\\
\label{gb}
\lim_{x\to n}\Gamma(-2x-m) &\approx \frac{(-1)^{1+2 n+m}}{2\Gamma(1+2 n+m) } \frac{1}{x-n} +O(1),\mbox{$\,\,m,\,n \in 0\cup \mathbb{Z}^+$,}\\
\label{gc}
\lim_{x\to n} \lim_{y\to m} \frac{\Gamma(-2x)}{\Gamma(-x-y)}& \approx (-1)^{n+m} \frac{\Gamma(1+n+m)}{2\Gamma(1+2n)} \left(1+\frac{x-n}{y-m} \right)\\ \nonumber
&+O((x-n),(y-m)),\mbox{ $\,\,m,\,n \in \mathbb{Z}^+$,} \\
\label{gd}
\lim_{x\to n} \lim_{y\to m} \frac{\Gamma(1-x)}{\Gamma(1-x-y)}& \approx (-1)^{m} \frac{\Gamma(n+m)}{\Gamma(n)} \left(1+\frac{y-m}{x-n} \right)\\ \nonumber
&+O((x-n),(y-m)),\mbox{ $\,\,m,\,n \in \mathbb{Z}^+$.}
\end{align} 
\end{subequations}

\renewcommand{\theequation}{B\arabic{equation}}
\setcounter{equation}{0}
\section{Alternative boundary conditions, Problem 1, Rule III}
\label{Z}

The fact that $Y^{(0)}$ fulfills the boundary condition (\ref{x2bc}) for $k=0$ can be verified through a simple substitution. In order to see that the boundary condition for $Y^{(k>0)}$ is also satisfied one needs to note that\footnote{We define $t \equiv -\sqrt{x_1}-\frac{1}{\sqrt{x_1}}$ and study the limit of $Y^{(k>0)}$ as $x_2 \to t$.} the overall multiplicative factor of (\ref{solkb}) decays as $(x_2-t)^{1+2k}$. Expanding the first term of the curly bracket (product of three factors) shows that this quantity grows as $(x_2-t)^{-k-1}$. Next, we move to the second term in the curly bracket, which is the term involving the hypergeometric function, and whose argument grows as $(x_2-t)^{-1}$. Taking into account that $k$ is a positive integer, we deduce that the $_2F_1\left(-1-2k,-k;-2k,x\right)$ is a terminating polynomial of degree $k$ and hence, we conclude that this term grows as  $(x_2-t)^{-k}$. Thus, the curly bracket grows as $(x_2-t)^{-k-1}$ while the overall coefficient decays as $(x_2-t)^{1+2k}$ and thus the whole equation decays as $(x_2-t)^{k}$. This implies that the boundary condition is fulfilled $\forall \,k>0$. 

We thus conclude that the solutions (\ref{solk}) could had been derived from the same differential equation (\ref{odep1a}) but with the alternative boundary conditions (\ref{x2bc}).

\renewcommand{\theequation}{C\arabic{equation}}
\setcounter{equation}{0}
\section{A second derivation of (\ref{ffinal2})}
\label{B}

The current appendix describes briefly a second derivation of equation (\ref{ffinal2}). Starting from (\ref{fb1}) and performing the summation using standard Gamma function identities we obtain
\begin{align} \label{fA1}
f_{r,w}  = \frac{(-1)^w 4^r  \Gamma(\frac{1}{2}+r)}{\sqrt{\pi}\Gamma(2+r)\Gamma(1+w)}  \times \frac{\Gamma(-2r)}{\Gamma(-r-w)}  \frac{\Gamma(1-w)}{ \Gamma(1-r-w)},\,\, (r,w)\neq(0,0).
\end{align}
Using (\ref{gc}) and (\ref{gd}) for respectively the last two Gamma function ratios, (\ref{fA1}) becomes
\begin{align} \label{fA2}
f_{r,w}  &= \frac{ 4^r  \Gamma(\frac{1}{2}+r)}{\sqrt{\pi}\Gamma(2+r)\Gamma(1+w)}   \frac{\Gamma(1+r+w)}{\Gamma(1+2r)}  
 \frac{\Gamma(r+w)}{ \Gamma(w)} + \lim_{x\to r,y\to w}\left(O\left(\frac{y-w}{x-r},\frac{x-r}{y-w}\right) \right)\\ \nonumber
&=\frac{\Gamma(r+w)\Gamma(1+r+w)}{\Gamma(1+r) \Gamma(2+r)\Gamma(w)\Gamma(1+w)}
\end{align}
where in arriving in the second equality we used the identity $\Gamma(2z) = \frac{1}{\sqrt{2\pi}}2^{2z-\frac{1}{2}}\Gamma(z)\Gamma(z+\frac{1}{2})$ and dropped the $O\left(\frac{y-w}{x-r},\frac{x-r}{y-w}\right) $ terms assuming some regulation argument. Equation (\ref{fA2}) then can be written in terms of factorials and when it is simplified, it becomes exactly identical to (\ref{ffinal2}).
It is notable that unless we took sub-leading corrections in the expansions as in (\ref{gc}) and (\ref{gd}), we would be off by a factor of 2 with respect to  (\ref{ffinal2}).
\renewcommand{\theequation}{D\arabic{equation}}
\setcounter{equation}{0}
\section{Finite summation identities involving a product of binomials and a rational function}
\label{C}

We show that for positive integers $r$ and $w$
\begin{align} \label{C1}
\sum_{n=0}^r (-1)^{r+n}\frac{ \Gamma(1+r+w+n)} { (1+n)\Gamma(1+w) \Gamma(1+r-n) \Gamma^2(1+n) } =\left(\frac{(r+w)!}{r! w!}\right)^2 \frac{w}{(r+1)(r+w)}.
\end{align}
The left hand side of (\ref{C1}) for integer parameters is
\begin{align} \label{C2}
\sum_{n=0}^r (-1)^{r+n} & \frac{ \Gamma(1+r+w+n)} { (1+n)\Gamma(1+w) \Gamma(1+r-n) \Gamma^2(1+n) } = (-1)^r \sum_{n=0}^r \frac{(-1)^{n}}{n+1} \frac{(r+w+n)!}{w!(r-n)! (n!)^2}\\ \nonumber
 & =(-1)^r  \binom{r+w}{w} \sum_{n=0}^r \frac{(-1)^{n}}{n+1} \binom{r}{n} \binom{r+w+n}{r+w}.
\end{align}
 
Next, we use the summation identity
\begin{align}\label{id1}
\sum_{n=0}^r (-1)^n \binom{r}{n} \binom{s+n}{m} =(-1)^r \binom{s}{m-r},
\end{align}
which implies
\begin{align}\label{id11}
\sum_{n=0}^r& \frac{(-1)^n}{1+n} \binom{r}{n} \binom{s+n}{m} = \sum_{n=0}^r \frac{(-1)^n}{1+r} \binom{r+1}{n+1} \binom{s+n}{m} \\ \nonumber
&= -\sum_{n=1}^{1+r} \frac{(-1)^{n}}{1+r} \binom{r+1}{n} \binom{s-1+n}{m} \\ \nonumber
&=-\frac{1}{1+r}  \Bigg \{ \sum_{n=0}^{1+r}\left( (-1)^{n} \binom{r+1}{n} \binom{s-1+n}{m} \right) - \binom{r+1}{0} \binom{s-1}{m} \Bigg \} \\ \nonumber
&=\frac{1}{1+r} \Bigg \{  (-1)^{r} \binom{s-1}{m-r-1} + \binom{s-1}{m}\Bigg \}. 
\end{align}
Inserting (\ref{id11}) in (\ref{C2}) with $s=m=r+w$ and noting that the second binomial in the most right hand side of (\ref{id11}) vanishes, we find
\begin{align} \label{CendProof}
\sum_{n=0}^r (-1)^{r+n} & \frac{ \Gamma(1+r+w+n)} { (1+n)\Gamma(1+w) \Gamma(1+r-n) \Gamma^2(1+n) } = 
\frac{1}{1+r} \binom{r+w}{w} \binom{r+w-1}{w-1} \\ \nonumber
& =\left(\frac{(r+w)!}{r! w!}\right)^2 \frac{w}{(r+1)(r+w)},
\end{align}
and this completes the proof.

\renewcommand{\theequation}{E\arabic{equation}}
\setcounter{equation}{0}
\section{Proving the normalization of the probability formula, Rule III}
\label{D}
Starting from (\ref{finallypk}) we prove that the normalization condition (\ref{normalp}) applies. Equation (\ref{finallypk}) yields

\begin{align}\label{normalpd1}
\sum_{k=0}^w p^{(k)}(r,w)& = \frac{rw!r!(r+w+1)}{(r+w)!} \sum_{k=0}^w \left \{ \frac{k! }{(r+k+1)!}\,\, \frac{(r+w-k-1)!}{(w-k)!} \right \} \\ \nonumber
&=   \frac{rw!r!(r+w+1)}{(r+w)!} \left \{  \frac{\Gamma(r+w)}{\Gamma(w+1) \Gamma(2+r)}\,  _3F_2\left(1,1,-w,2+r,1-r-w;1 \right)  \right \} 
\end{align}  
where the summation is, up to overall Gamma function factors, a hypergeometric function $_3F_2$. The next step is to use the identity \cite{wolfram}
\begin{align}\label{did}
_3F_2(a,b,-n,d,a+b-d-n+1;1) = \frac{(d-a)_n (d-b)_n }{(d)_n (d-a-b)_n}, \, n \in \mathbb{N}
\end{align}
for $a=b=1$, $n=w$, $d=2+r$, which yield $a+b-d-n+1 = 1-r-w$, and where $x_n \equiv \Gamma(x+n)/\Gamma(x)$ is the Pochhammer's Symbol. This means that the identity is applicable for equation (\ref{normalpd1}), which after simplifying yields
\begin{align}\label{normalpd2}
\sum_{k=0}^w p^{(k)}(r,w) =   \frac{r r!(1+r+w)}{(1+r)!(r+w)} \,  \frac{\left((1+r)_w\right)^2  }{(2+r)_w (r)_w} = 1
\end{align}  
completing the proof.

\renewcommand{\theequation}{F\arabic{equation}}
\setcounter{equation}{0}
\section{Solution to Problem 1, Rule III: a combinatorial-inductive derivation}
\label{Gminus}
In this appendix, a combinatorial-inductive derivation of (\ref{finallypk}), which is shorter and more intuitive, is given. On the other hand, this derivation is less systematic and it exclusively refers to the particular problem (Problem 1). In particular,  it may not be generalized and adapted to other similar problems. An analogous derivation, as the one presented here using a similar logic, exists for Problem 1, Rule IV and its presentation is therefore omitted.

We consider sequences containing $r$ red ball removals out of $r$ and $w-k$ white ball removals out of $w$, ending with a red ball removal. Each removal occurs with probability $r^2/(r+w)^2$ for a red and with probability $w(2r+w)/(r+w)^2$ for a white ball respectively, where $r$ and $w$ are the number of the red and of the white balls at the time of the given removal. Defining with 
\begin{align}\label{ri}
r_i& \equiv \mbox{The number of the red balls present during the $i$-th white ball removal,}\\ \notag
i&=1,\,2,\,...,\, w-k,\mbox{ with $r\geq r_1\geq r_2\geq r_3\geq...\geq r_{w-k}\geq 1$,}
\end{align}
the probability of a given sequence of removals has the generic form
\begin{align}\label{seq}
\left(\frac{r^2}{(r+w)^2}\frac{(r-1)^2}{((r-1)+w)^2}...\frac{(r_1+1)^2}{((r_1+1)+w)^2}\right) \frac{w(2r_1+w)}{(r_1+w)^2} \times\\ \notag
 \left(\frac{r_1^2}{(r_1+(w-1))^2}\frac{(r_1-1)^2}{((r_1-1)+(w-1))^2}...\frac{(r_2+1)^2}{((r_2+1)+(w-1))^2}\right)  \frac{(w-1)(2r_2+(w-1))}{(r_2+(w-1))^2} \times \\ \notag
 \left(\frac{r_2^2}{(r_2+(w-2))^2}\frac{(r_2-1)^2}{((r_2-1)+(w-2))^2}...\frac{(r_3+1)^2}{((r_3+1)+(w-2))^2}\right)  \frac{(w-2)(2r_3+(w-2))}{(r_3+(w-2))^2}\times \\ \notag
... \frac{(r_{w-k}+1)^2}{((r_{w-k}+1)+(k+1))^2} \frac{(k+1)(2r_{w-k}+(k+1))}{(r_{w-k}+(k+1))^2}\times  \\ \notag
\frac{r_{w-k}^2}{(r_{w-k}+k)^2}  \frac{(r_{w-k}-1)^2}{((r_{w-k}-1)+k)^2} ....\frac{1}{(1+k)^2}.
\end{align}
It could happen that sometimes the sub-product of terms between successive white ball removals, i.e. the terms between $\frac{(w-i+1)(2r_i+(w-i+1))}{(r_i+(w-i+1))^2}$ and $\frac{(w-i)(2r_{i+1}+(w-i))}{(r_{i+1}+(w-i))^2}$, $i=1,\,2,\,...,\,w-k-1$, including the first sub-product of terms until the $\frac{w(2r_1+w)}{(r_1+w)^2}$ factor, could collapse to unity. Equation (\ref{seq}) can be written as a product of two factors according to
\begin{align}\label{seq1}
\Bigg[\left(\frac{r^2}{(r+w)^2}\frac{(r-1)^2}{((r-1)+w)^2}...\frac{(r_1+1)^2}{((r_1+1)+w)^2}\right) \frac{w}{(r_1+w)^2} \times\\ \notag
 \left(\frac{r_1^2}{(r_1+(w-1))^2}\frac{(r_1-1)^2}{((r_1-1)+(w-1))^2}...\frac{(r_2+1)^2}{((r_2+1)+(w-1))^2}\right)  \frac{w-1}{(r_2+(w-1))^2} \times \\ \notag
 \left(\frac{r_2^2}{(r_2+(w-2))^2}\frac{(r_2-1)^2}{((r_2-1)+(w-2))^2}...\frac{(r_3+1)^2}{((r_3+1)+(w-2))^2}\right)  \frac{w-2}{(r_3+(w-2))^2}\times \\ \notag
... \frac{(r_{w-k}+1)^2}{((r_{w-k}+1)+(k+1))^2} \frac{k+1}{(r_{w-k}+(k+1))^2}\times  \\ \notag
\frac{r_{w-k}^2}{(r_{w-k}+k)^2}  \frac{(r_{w-k}-1)^2}{((r_{w-k}-1)+k)^2} ....\frac{1}{(1+k)^2} \Bigg] \times \prod_{i=1}^{w-k} \left(2r_i+(w+1-i) \right)  \\ \notag
= \left[ \frac{(r!)^2 \prod_{i=1}^{r+w-k} (w+1-i)}{\left(\prod_{i=1}^{r+w-k} (r+w+1-i)\right)^2} \right] \times \prod_{i=1}^{w-k} \left(2r_i+(w+1-i) \right)   \\ \notag
 = \frac{(r!)^2 k! w!}{((r+w)!)^2} \times\prod_{i=1}^{w-k} \left(2r_i+(w+1-i) \right)
\end{align}
where the first factor is the term in the square bracket, which simplifies to $ \frac{(r!)^2 k! w!}{((r+w)!)^2}$. Interestingly, all the $r_i$ dependance drops out from this first term and hence, this factor becomes a multiplicative overall factor to the sequence (\ref{seq}).

Using (\ref{ri}) and (\ref{seq1}), the required probability is then given by
\begin{align}\label{psum}
p_{III}^{(k)}(r,w)& =  \frac{(r!)^2 k! w!}{((r+w)!)^2} \sum_{r\geq r_1\geq r_2\geq...\geq r_{w-k}\geq 1} (2r_1+w)(2r_2+(w-1))...(2r_{w-k}+(k+1))\\ \notag
&=  \frac{(r!)^2 k! w!}{((r+w)!)^2} \sum_{r\geq r_1\geq r_2\geq...\geq r_{s}\geq 1} (2r_1+k+s)(2r_2+k+(s-1))...(2r_{s}+(k+1))\\ \notag
&=  \frac{(r!)^2 k! w!}{((r+w)!)^2} g^{(k)}(r,s)
\end{align}
 where we have set $w=k+s$ and where we have defined
 \begin{align}\label{gk}
 g^{(k)}(r,s) \equiv \sum_{r\geq r_1\geq r_2\geq...\geq r_{s}\geq 1} (2r_1+k+s)(2r_2+k+(s-1))...(2r_{s}+(k+1)).
 \end{align}
Starting from the definition of $ g^{(k)}(r,s)$, equation (\ref{gk}), it is straightforward to show
 \begin{align}\label{bcg}
g^{(k)}(1,s) = \frac{(k+s+2)!}{(k+2)!},\,\,\, \,\,g^{(k)}(r,1) = r(r+k+2) .
 \end{align}
So the problem reduces in computing $g^{(k)}(r,s=w-k)$ from (\ref{gk}) subject to the boundary cases (\ref{bcg}). 
We proceed by taking the difference $g^{(k)}(r,s)-g^{(k)}(r-1,s)$ from where we find
 \begin{align}
g^{(k)}(r,s)-g^{(k)}(r-1,s) = \hspace{-0.5cm} \sum_{r= r_1\geq r_2\geq...\geq r_{s}\geq 1} \hspace{-0.5cm}(2r_1+k+s)(2r_2+k+(s-1))...(2r_{s}+(k+1))\\ \notag
\hspace{-0.5cm}= (2r+k+s) \hspace{-0.5cm} \sum_{r\geq r_2\geq...\geq r_{s}\geq 1} \hspace{-0.5cm}(2r_2+k+(s-1))...(2r_{s}+(k+1))=  (2r+k+s) g^{(k)}(r,s-1)
 \end{align}
 where in the second equality we have assumed a renaming of the dummy $r_i$ index according to $r_i \rightarrow r_{i-1}$, $i=2,\,3\,...,\,s$. Thus we have arrived at the recursion\footnote{If we had used the alternative ansatz $p_{III}^{(k)}(r,w) =  \frac{(r!)^2 k! w!}{((r+w)!)^2} f^{(k)}(r,w)$ instated of the ansatz of (\ref{ans2}), equation (\ref{rec3}) would had yielded the much simpler difference equation $f^{(k)}(r,w)=f^{(k)}(r-1,w)+(2r+w)f^{(k)}(r,w-1)$ rather than the recursion provided by equation (\ref{rec4}).} 
 \begin{align}\label{grec}
 g^{(k)}(r,s)-g^{(k)}(r-1,s) =  (2r+k+s) g^{(k)}(r,s-1)
 \end{align}
subject to equations (\ref{bcg}). In the remaining part we will only sketch the derivation rather than provide the details. While we could proceed by induction, we follow a more straightforward approach.

Using (\ref{bcg}) and (\ref{grec}) for $s=2$ it is deduced that
\begin{equation}
g^{(k)}(2,2)=    3 (k+4)(k+5),\,\,
g^{(k)}(3,2)=  6 (k+5)(k+6),\,\,
g^{(k)}(4,2)=   10 (k+6)(k+7),\,\,
\end{equation}
from where it is guessed that
\begin{align}\label{gs2}
g^{(k)}(r,2) = \frac{r(r+1)}{2} (k+r+2)(k+r+3).
\end{align}
Equation (\ref{gs2}) can then be substituted in (\ref{grec}), and with the help of (\ref{bcg}), one can verify that (\ref{grec}) is indeed fulfilled. The next step is to work for $s=3$ using (\ref{bcg}) and (\ref{gs2}) to obtain
\begin{align}
g^{(k)}(2,3)&=    4 (k+4)(k+5) (k+6),\,\,\,\,\,
g^{(k)}(3,3)=  10 (k+5)(k+6)(k+7),\,\,\\ \notag
g^{(k)}(4,3)&=   20 (k+6)(k+7)(k+8),\,\,
g^{(k)}(5,3)=   35 (k+7)(k+8)(k+9)
\end{align}
from where it is guessed that
\begin{align}\label{gs3}
g^{(k)}(r,3) = \frac{r(r+1)(r+2)}{3!} (k+r+2)(k+r+3)(k+r+4).
\end{align}
Equation (\ref{gs3}) can then be verified through a substitution in (\ref{grec}) with the help of (\ref{gs2}). 

Observing then equations (\ref{gs2}) and (\ref{gs3}) we make the general ansatz solution
\begin{align}\label{gsr}
g^{(k)}(r,s) = \binom{r+s-1}{s} \frac{(r+s+k-1)!}{(r+k+1)!},
\end{align}
which we verify by direct substitution in the recursion (\ref{grec}). 

Finally, using (\ref{psum}) and (\ref{gsr}) we obtain
\begin{align}\label{p3com}
p_{III}^{(k)}(r,w)& =  \frac{(r!)^2 k! w!}{((r+w)!)^2} g^{(k)}(r,s=w-k) = \frac{k! r!(r+w+1)}{(r+k+1)!}\,\, \frac{rw!(r+w-k-1)!}{(r+w)!(w-k)!}
\end{align}
reproducing equation (\ref{finallypk}), and this completes the alternative derivation of $p_{III}^{(k)}(r,w)$.

\underline {Remark:} We could had started from  (\ref{rec4}) and subsequently worked along the same lines as between equations (\ref{grec})-(\ref{p3com}). Conversely, starting from (\ref{grec}), we could had worked along the lines between equation (\ref{rec4}) and equation (\ref{finallypk}); that is using the PDE approach. Either approach has its advantages: using the method of this appendix we reached the result faster. On the other hand, the approach here has involved educated guessing, which (in principle) has not been the case during the initial derivation presented by the main body of this paper. That approach has been purely systematic and it may, in principle, be generalized and be applied to other cases, as that approach is to a large extend problem independent. 

\renewcommand{\theequation}{G\arabic{equation}}
\setcounter{equation}{0}
\section{Problem 4: A challenging problem from the literature}
\label{E}

The problem presented in \cite{OaPe} is as follows.

\vspace{0.2cm}
\hspace{-0.7cm}{\bf Problem 4}
\vspace{0.2cm}

\vspace{0.2in}

\hspace{-0.7cm}
{\it A bag contains $m > 0$ black balls and $n > 0$ white balls. A sequence of balls from the bag is discarded in the following manner:

\vspace{0.5in}

(i) A ball is chosen at random and discarded. (ii) Another ball is chosen at random from the remainder. If its color is different from the last it is replaced in the bag and the process repeated from the beginning (i.e. (i)). If the second ball is the same color as the first it is discarded and we proceed from (ii). Thus the balls are sampled and discarded until a change in color occurs, at which point the last ball is replaced and the process starts afresh.

The question is: what is the probability that the final ball should be black?
}

\vspace{0.2cm}
For the record, the authors show that the requested probability is as simple as $1/2$.

\renewcommand{\theequation}{H\arabic{equation}}
\setcounter{equation}{0}
\section{Python simulation code and some results}
\label{F}
The python codes that simulate the removals of Problem 1 are provided below. 

The first block of the code that follows (Problem 1, Rule III) imports several necessary Python packages. The ``\texttt{newListAfterRemoval}" function provides the updated content of the box once a ball is removed. The function ``\texttt{redWhite\_Problem1\_Rule\_IIII}($r,w$,epochs)" refers to Problem 1, Rule III, and it takes care of everything else. The arguments $r$ and $w$ refer to the initial red ($r$) and the initial white ($w$) balls respectively. The variable ``epochs" is a positive integer and defines the number of times that the experiment will be repeated. In particular, this function provides the simulated probabilities (Rule III) contrasting them with the theoretical ones. This part of the code is supplemented by explanations (see comments in green fonds).

Likewise, the function ``\texttt{redWhite\_Problem1\_Rule\_IV}($r,w$,epochs)" refers to Problem 1, Rule IV, and it takes care of everything else. The arguments in the function are analogous to those of Rule III. In particular, this function provides the simulated probabilities (Rule IV) contrasting them with the theoretical ones. This part of the code is not supplemented by explanations as these are analogous to those of ``\texttt{redWhite\_Problem1\_Rule\_IIII}($r,w$,epochs)".


\includepdf[pages=1,pagecommand={},offset=1.9cm -4cm]{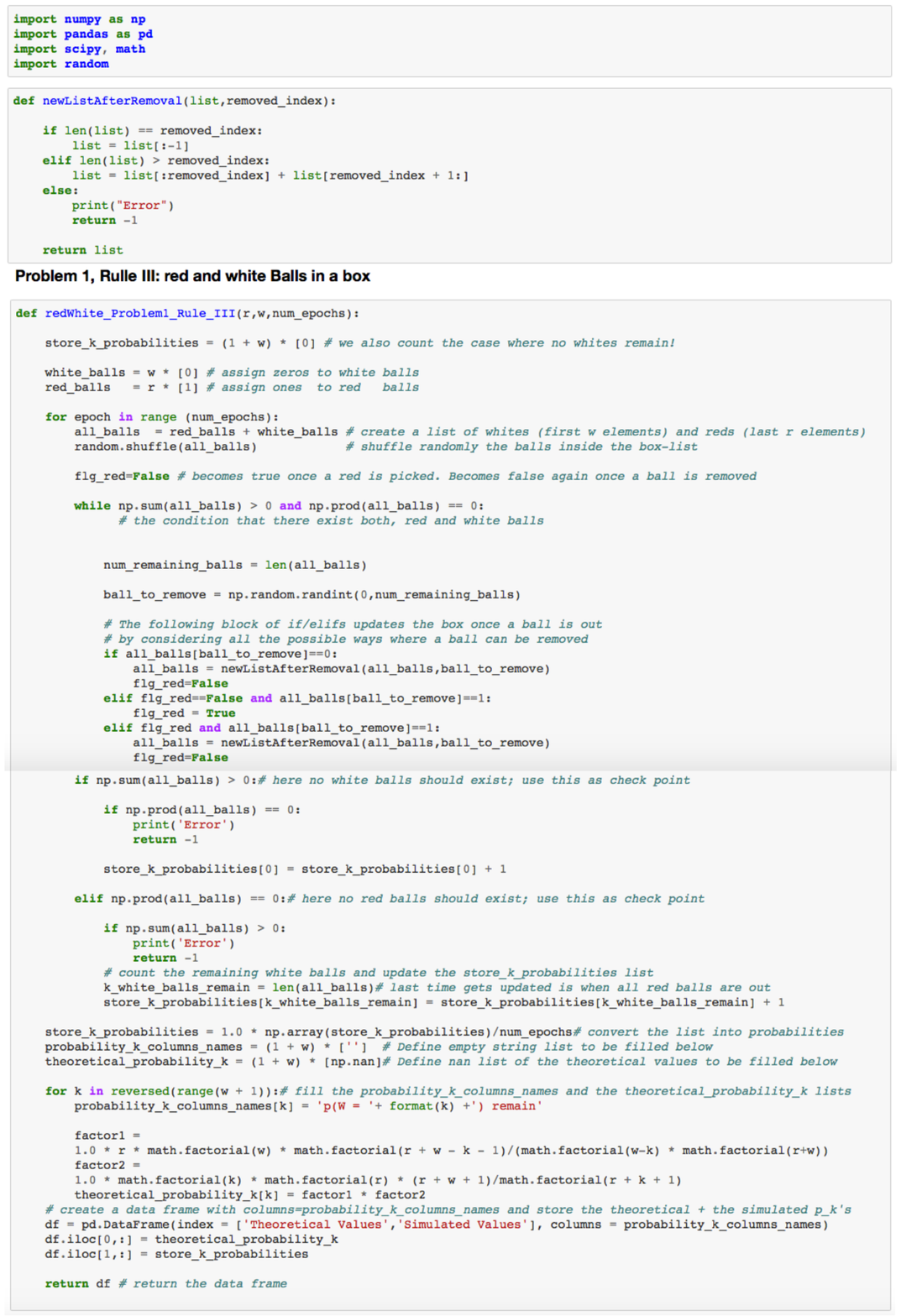}

Below we present the simulated and the theoretical results for a few cases for  Problem 1, Rule III.
\vspace{-0in}
\begin{figure}[!th]
\centering
\hspace{-1in}\includegraphics*[scale=0.95,clip=true]{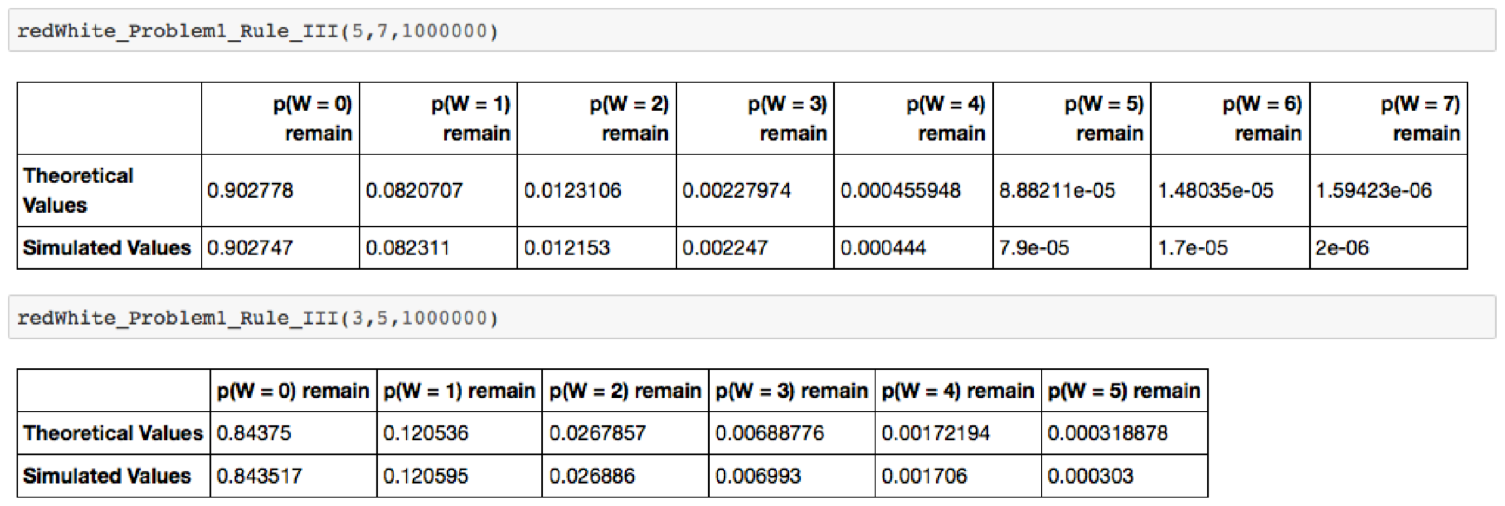} 
\vspace{-7in}
\caption{Problem 1, Rule III. The simulated Vs the theoretical results for 1 000 000 simulations in two cases: Upper panel refers to $r=5$ and $w=7$. Lower panel refers to $r=3$ and $w=5$. It is observed that the probability distribution is maximized when $k=0$ in accordance with (\ref{maxpk3b}). The simulated values are very close to the theoretical ones especially for smaller values of $k$ (smaller number of remaining white balls).}
\label{resultsIII}
\end{figure}

Next, we present the code and the results for Problem 1, Rule IV.

\includepdf[pages=1,pagecommand={},offset=1.9cm -4cm]{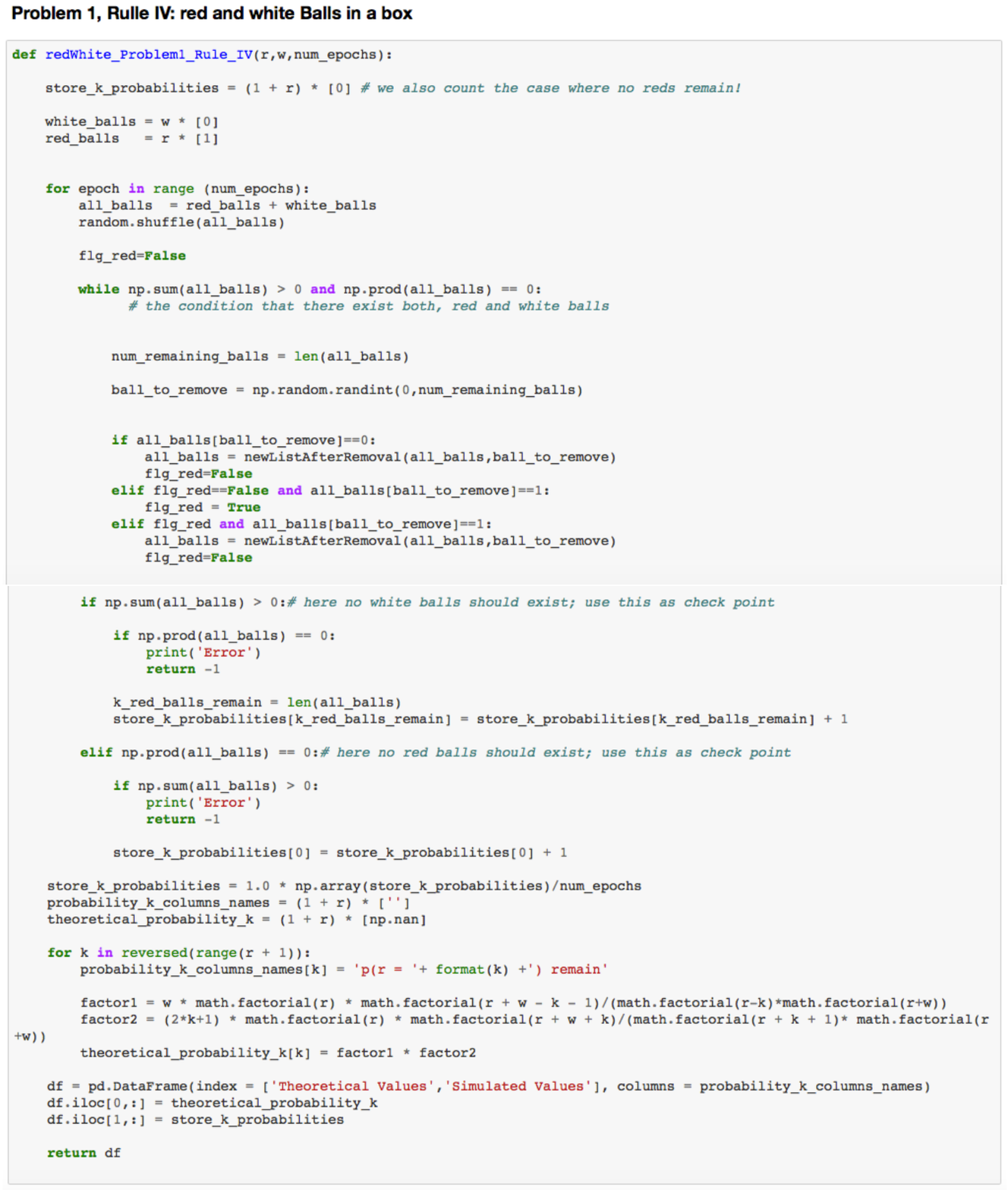}

Below we present the simulated and the theoretical results for a few cases for  Problem 1, Rule IV.
\vspace{-1.0in}

\begin{figure}[!th]
\centering
\includegraphics*[scale=0.95,angle=0,clip=true]{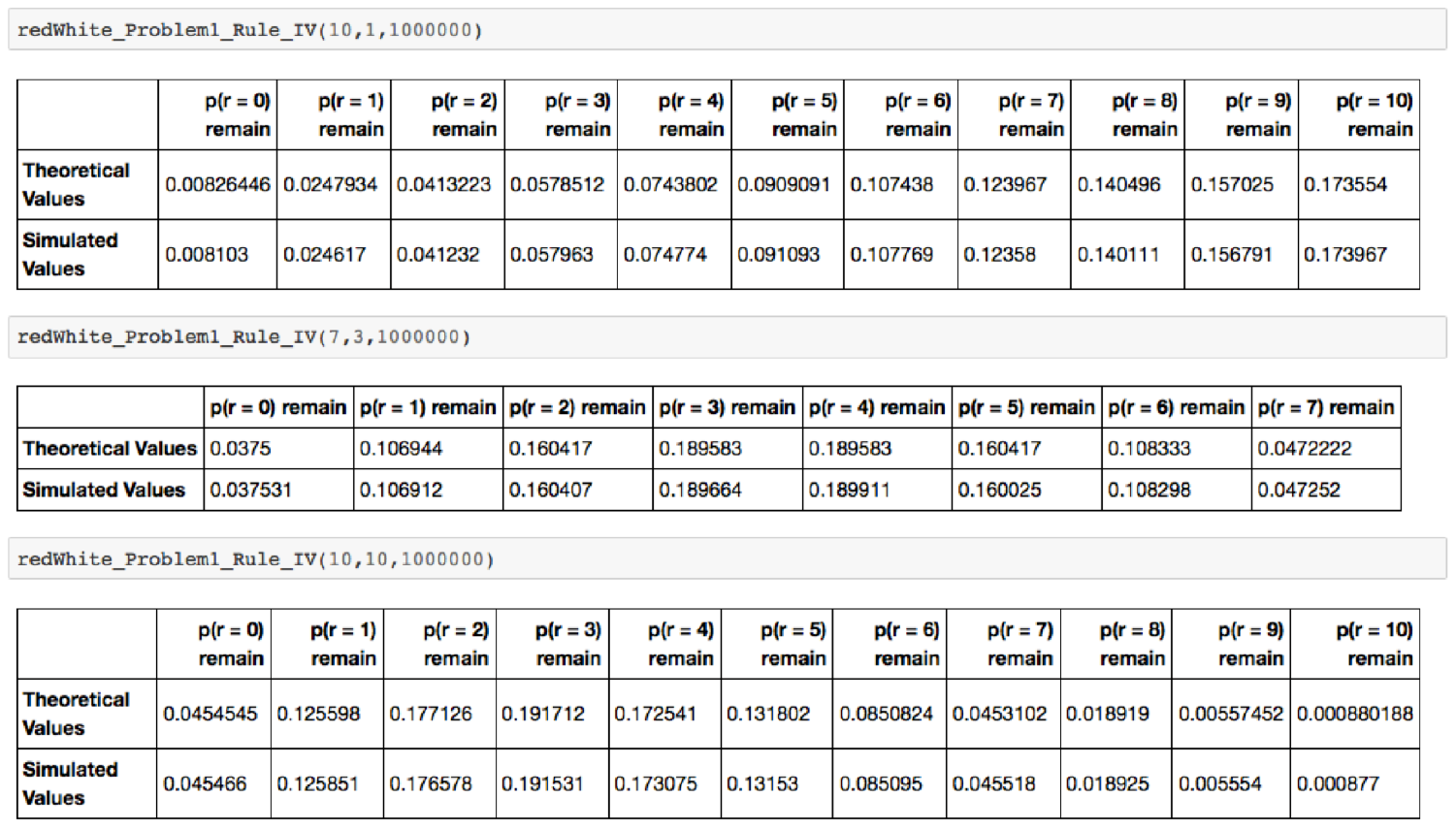} 
\vspace{-6in}
\caption{Problem 1, Rule IV. The simulated Vs the theoretical results for 1 000 000 simulations in three cases: 
Upper panel refers to $r=10$ and $w=1$. This case corresponds and confirms equation (\ref{maxpk4c}) with a unique maximum at $k=r=10$.
Middle panel refers to  $r=7$ and $w=3$. This case corresponds to $k_0=4 \in \mathbb{Z}_+ $ (see (\ref{k0})) and confirms equation (\ref{maxpk4d}) with degenerate maxima at $k=3$ and $k=4$.
Lower panel refers to $r=10$ and $w=10$. This case corresponds to $k_0=2\sqrt{\frac{55}{19}} \approx 3.4 \notin \mathbb{Z}_+ $ and confirms equation (\ref{maxpk4e}) with a unique maximum at $k=3$.
It is observed that the simulated values are very close to the theoretical ones.}
\label{resultsIV}
\end{figure}



\addcontentsline{toc}{section}{References}
\providecommand{\href}[2]{#2}\begingroup\raggedright\endgroup

\end{document}